\documentclass{amsart}
\usepackage{amssymb,stmaryrd}
\usepackage{amsfonts}
\usepackage{amstext}
\usepackage{algorithmic}
\usepackage{algorithm}
\usepackage{graphicx}
\usepackage{epstopdf}
\usepackage[all]{xy}
\usepackage{MnSymbol}
\numberwithin{equation}{section}
\newtheorem{theorem}{Theorem}[section]

\newtheorem{proposition}[theorem]{Proposition}

\theoremstyle{definition}

\theoremstyle{remark}

\newcommand{\R}{{\mathbb{R}}}

\newcommand{\C}{{\mathbb{C}}}
\newcommand{\Z}{{\mathbb{Z}}}

\newcommand{\F}{{\mathbb{F}}}
\newcommand{\N}{{\mathbb{N}}}

\renewcommand{\>}{{\rangle}}

\newcommand{\wedgeq}{{\wedge\kern-5pt\cdot}}

\renewcommand{\ker}{{\rm{ker}}}

\newcommand{\tens}{\otimes}

\newcommand{\id}{{\rm id}}

\newcommand{\extd}{{\rm d}}
\newcommand{\del}{{\partial}}

\begin{document}
\title{Classification of digital affine noncommutative geometries}
\author{Shahn Majid \& Anna Pacho{\l }}
\address{School of Mathematics, Mile End Rd, London E1 4NS, UK}
\thanks{
This project has received funding from the European Union's Horizon 2020
research and innovation programme under the Marie Sk{\l }odowska-Curie grant
agreement No 660061}

\email{s.majid@qmul.ac.uk, a.pachol@qmul.ac.uk}

\begin{abstract}
It is known that connected translation invariant $n$-dimensional
noncommutative differentials $\extd x^i$ on the algebra $k[x^1,\cdots,x^n]$ of
polynomials in $n$-variables over a field $k$ are classified by commutative
algebras $V$ on the vector space spanned by the coordinates. This data
also applies to construct differentials on the Heisenberg algebra
`spacetime' with relations $[x^\mu,x^\nu]=\lambda\Theta^{\mu\nu}$ where $
\Theta$ is an antisymmetric matrix as well as to Lie algebras with pre-Lie
algebra structures. We specialise the general theory to the field $k={\ 
\mathbb{F}}_2$ of two elements, in which case 
translation invariant metrics (i.e. with constant coefficients) are equivalent to making $V$ a Frobenius algebras. 
We classify all of these and their quantum Levi-Civita
bimodule connections for $n=2,3$, with partial results for $n=4$. For $n=2$ we
find 3 inequivalent differential structures admitting 1,2 and 3 invariant
metrics respectively. For $n=3$ we find 6 differential structures admitting $0,1,2,3,4,7$
invariant metrics respectively. We give some examples for $n=4$ and general $n$. Surprisingly, not all our geometries for $n\ge 2$ have zero quantum Riemann curvature. Quantum gravity is
normally seen as a weighted `sum' over all possible metrics but our results
are a step towards a deeper approach in which we must also `sum' over
differential structures.  Over ${\mathbb{F}}_2$ we
construct some of our algebras and associated structures by digital gates, opening
up the possibility of `digital geometry'. \end{abstract}

\maketitle

\section{Introduction}

A standard technique in physics and engineering is to replace geometric backgrounds by discrete approximations
such as a lattice or graph, thereby rendering systems more calculable. In recent years it has become clear that
this can be handled by noncommutative geometry not because the `coordinate algebras' $A$ are noncommutative (they remain
commutative) but because differentials and functions do not commute, see \cite{Ma:gra} and references therein. The formalism of noncommutative differential geometry does not require functions and differentials to commute, so is more general even when the
algebra is classical. In the present work we use such noncommutative differential geometry to 
explore a different and complementary kind of `discretisaton scheme' in which the field $\C$ or $\R$ that we
work over is replaced by the field $\F_2$ of two elements $0,1$ and which we call {\em digital geometry}. 

We use a `bottom up' constructive approach to noncommutative differential geometry that grew in the 1990s out of (but not limited to) the differential geometry of quantum groups, rather than one of the powerful operator algebra approach to noncommtutative geometry as in \cite{Con}. This is more explicit
(albeit mathematically less deep) and has the merit that one can work over any field $k$. Often characteristic $2$ (which includes $\F_2$) is excluded for simplicity so one must be a little careful (notably tensors cannot be decomposed into symmetric and antisymmetric parts) but
most of the theory including differential forms (as differential graded algebras $\Omega(A)$), vector bundles, principal bundles, connections
and Riemannian metrics work over any field. We refer to our LTCC lectures \cite{Ma:ltcc} for a recent introduction. A small part of the formalism is recapped in Section~\ref{secformalism} along
with a recent classification theorem \cite{MaTao} for translation invariant differentials on Hopf algebras with linear (additive) coproduct, 
which will be our starting point. To keep a lid on the classification problem we insist that our metrics are invertible, which is known \cite{BegMa2}
to require that the metric is central (commutes with functions) and we assume that our connections $\nabla:\Omega^1\to \Omega^1\tens_A\Omega^1$ are {\em bimodule connections} \cite{Mou,DVM}. This means that their right handed derivation rule
is expressed in a `generalised braiding' $\sigma:\Omega^1\tens_A\Omega^1\to \Omega^1\tens_A\Omega^1$ and we require this to be invertible. The `quantum groups' approach to noncommutative differential geometry was particularly developped using bimodule connections in recent works such as \cite{BegMa1,BegMa2,Ma:gra,BasMa}. 

The present paper follows on from \cite{BasMa} where we studied the de Rham cohomology of $\F_2[x]$ (polynomials in one variable) with noncommutative differential structures, which turned out to be surprisingly rich. This led to  nice family of `finite' geometries over $\F_2$ as finite dimensional commutative Hopf algebras $A_d$ for every $d\in\N$ (and over $\F_p$ for any prime $p$). By contrast, we will now be interested in affine or `flat space' $A=\F_2[x^1,\cdots,x^n]$ but it turns out that the classification of its differential structures of dimension $n$ already amounts to the classification of finite geometries in the form of $n$-dimensional commutative algebras $V$ over $\F_2$, so our results now include the classification of all of these up to dimension $n=4$ (we find that there are 16 of these up to isomorphism if we ask for them to be unital, see Section~\ref{secn4}), along with more complete results for $n=2,3$ for metrics and connections on $A$  for each of the respectively $3$ and $6$ possible choices of $V$ in these dimensions, see Sections~\ref{secn2}, \ref{secn3}. These bimodule noncommutative geometries are explored under the restriction that the metric and Christoffel symbol coefficients are constants in keeping with our view of $A=\F_2[x^1,\cdots,x^n]$ as `flat space', i.e. we are looking at  translation invariant geometries over $\F_2$. It is interesting that for $n\ge 3$ some of the possible geometries nevertheless have quantum curvature $R_\nabla\ne 0$, which we regard as a purely quantum phenomenon.

We envisage many applications throughout mathematical physics and engineering wherever classical differential geometry plays a role. It is not our goal to develop these here but we conclude with an extended discussion in Section~\ref{secdis} of some that we have in mind. Our own motivation for noncommutative geometry has come from quantum gravity in which the proposal of {\em quantum spacetime} and concrete models \cite{Ma:pla,MaRue,AmeMa} emerged out of quantum groups (and was the origin of one of the two main classes of quantum groups, namely the bicrossroduct ones). In this context one could in principle `sum over all geometries' so our classification is a peek into a restricted part of this. More generally our classification is a tool for model building and one can explore each of our geometries much further, for example solving wave equations. Clearly we would like to go further and explore all geometries not just the translation invariant ones in the present work. Also emebedded in our above explanation and surprisingly forced on us by  translation invariance of the differentials $\extd x^i$ is the set up of classical and quantum field theory in which we work with the space of functions on a linear space $V$ which is itself the space of functions on an underlying geometry. This suggests a different envisaged application in which spacetime would be the coordinate algebra $V$ and $A=k[x_1,\cdots,x_n]$ or more abstractly $k[V]$ would be the algebra of functionals {\em on} $V$ as the vector space of functions. Differentials also automatically extend to $A$ the Heisenberg algebra, so the first steps of quantum field theory also arise out of the natural possibilities for the noncommutative geometry of affine spaces. In this case our spacetime geometry is built on differentials and Riemannian structures on $\Omega(V)$ as in \cite{BasMa} and as will be classified in a sequel \cite{MaPac} in preparation. The discussion ends with a translation of algebra over $\F_2$ into digital electronics, thereby justifying our terminology and opening up a new front of applications in which geometric ideas can be translated into electronics.

We made extensive use of the numerical package R to enumerate all possible values of our structure constants, preceded and followed by symbolic calculations on Mathematica.

\section{Calculi on $k[x^1,\cdots,x^n]$ and Heisenberg algebras}
\label{secformalism}

If $A$ is a possibly noncommutative `coordinate' algebra, by differential
calculus on $A$ we mean an $A$-bimodule $\Omega^1$ and a map $\mathrm{d}
:A\to \Omega^1$ obeying the Leibniz rule $\mathrm{d}(ab)=(\mathrm{d} a)b+a 
\mathrm{d} b$ with the map $A\otimes A\to \Omega^1$ given by $a\otimes
b\mapsto a\mathrm{d} b$ surjective. Here a bimodule means we can
associatively multiply such 1-forms by elements of $A$ from the left and the
right. The calculus is called connected if $\ker\mathrm{d}=k.1$ where we
work over the field $k$. If $A$ is a Hopf algebra or `quantum group' the
coproduct expresses `group translation' and there is a standard notion of
the differential calculus being left and right covariant under this. We
refer to \cite{Ma:ltcc} for an introduction.

We build on the Majid-Tao theorem \cite{MaTao} which states that
connected translation invariant differential structures of classical
dimension on `quantum spaces' consisting of enveloping algebras $U(\mathfrak{g})$ where 
$\mathfrak{g}$ is a Lie algebra, are classified by pre-Lie structures $\circ 
$ on $\mathfrak{g}$. A pre-Lie algebra structure is a `product' $\circ $ on $
\mathfrak{g}$ such that 
\begin{equation}
v\circ w-w\circ v=[v,w]  \label{prelie1}
\end{equation}
(i.e. we recover the given Lie bracket) and 
\begin{equation}
(v\circ w)\circ z=(v\circ z)\circ w+v\circ (w\circ z-z\circ w),\quad \forall
v,w,z\in {g}.  \label{prelie2}
\end{equation}
The differential calculus has generators $\mathrm{d}x^{\mu }$ where $
\{x^{\mu }\}$ is a basis of $\mathfrak{g}$ and bimodule relations 
\begin{equation*}
\lbrack \mathrm{d}x^{\mu },x^{\nu }]=\lambda \mathrm{d}(x^{\mu }\circ x^{\nu
}).
\end{equation*}
Clearly the Jacobi identity 
\begin{equation*}
\lbrack \lbrack \mathrm{d}x^{\mu },x^{\nu }],x^{\rho }]=[[\mathrm{d}x^{\mu
},x^{\rho }],x^{\nu }]+[\mathrm{d}x^{\mu },[x^{\nu },x^{\rho }]]
\end{equation*}
for the bimodule relations translates immediately in this context to (\ref
{prelie2}). The Leibniz rule $\mathrm{d}[x^{\mu },x^{\nu }]=[\mathrm{d}
x^{\mu },x^{\nu }]+[x^{\mu },\mathrm{d}x^{\nu }]$ is the other part (\ref
{prelie1}). If the pre-Lie algebra is unital with identity element $e$ then
clearly the calculus is inner in the sense of existence of element $\theta
\in \Omega ^{1}$ such that $\mathrm{d}=[\theta ,]$ on $A$ (and on forms if
we use graded commutator), with $\theta =\lambda ^{-1}\mathrm{d}e$. Note
that the calculus could be inner in some other way with $\theta $ not the
differential of an element of the pre-Lie algebra. Isomorphisms of the
pre-Lie algebra are induced by linear coordinate transformations that do not
change the differential structure.

In the commutative case of $k[x^{1},\cdots ,x^{n}]$ regarded as the
enveloping algebra of an Abelian Lie algebra, we need $\circ $ commutative
and in this case (\ref{prelie2}) says that $\circ $ is associative, so the
data is that of an $n$-dimensional commutative algebra. Since $\mathrm{d}1=0$
and $1$ is central, a quick look at the proof above tells us that this works
just as well for the Heisenberg algebra regarded as noncommutative space, 
\begin{equation*}
\lbrack x^{\mu },x^{\nu }]=\lambda \Theta ^{\mu \nu }
\end{equation*}
which includes the commutative case with $\Theta =0$ (this can also be seen
as a Lie algebra with a central generator on the right hand side to which we
apply the pre-Lie theory and then set the central generator to 1).

There is in fact no need for $\mathfrak{g}$ to be finite dimensional. It can
be an infinite dimensional vector space $V$ with an antisymmetric bilinear
form $\Theta :V\times V\rightarrow k$ and the data for a calculus of the
above form on the associated algebra with relations $[v,w]=\lambda \Theta
(v,w)$ is precisely products $\circ :V\times V\rightarrow V$ making $
(V,\circ )$ an associative commutative algebra. In the unital case this will
be inner as before. An example is $V=C^{\infty }(M)$ on a manifold $M$ in
which case the above is a canonical noncommutative differential calculus or
`noncommutative variational calculus' on the space of functionals on $V$, or
more precisely on the symmetric algebra $S(V)$ or its Heisenberg `quantum
field theory' version.

Next we consider quantum metrics. In the constructive approach to
noncommutative geometry this means a nondegenerate element $g\in \Omega ^{1}\otimes
_{A}\Omega ^{1}$ which commutes with elements of $A$. The latter is known 
\cite{BegMa2} to be necessary for the existence of a bimodule inner
product $(\ ,\ ):\Omega ^{1}\otimes _{A}\Omega ^{1}\rightarrow A$ inverse to 
${g}$. We can also construct the latter directly. In our case a bimodule
inner product has the form $(\mathrm{d}v,\mathrm{d}w)=B(v,w)$ for some
bilinear map $V\times V\rightarrow A$ obeying 
\begin{equation}
B(v\circ w,z)+B(v,w\circ z)=\lambda ^{-1}[B(v,z),w]  \label{metric}
\end{equation}
for all $v,w,z\in V$. Here we require 
\begin{equation*}
((\mathrm{d}v)w,\mathrm{d}z)=wB(v,z)+\lambda B(v\circ w,z)=B(v,z)w-\lambda
B(v,z\circ w)=(\mathrm{d}v,w\mathrm{d}z)
\end{equation*}
using the commutation relations of differentials and functions, which in the
middle is the condition stated. For the inner product to be `real' we need $
\overline{B(v,w)}=B(w^{\ast },v^{\ast })$ which in a self-adjoint basis with 
$B$ symmetric requires its coefficients to be real.

However, if $\circ$ has an identity $e$ (so that the calculus is inner by a
coordinate differential) and the field has characteristic not $2$ then there
is no such map other than $B=0$ at the algebra level. So see this we set $w=e
$ so that $B(v,z)={\frac{1}{2\lambda}}[B(v,z),e]$ from our condition, which
has no solution at an algebraic level since the second expression has
strictly lower degree when $B(v,z)$ is written in a standard normal-ordered
form. There could still be non-algebraic examples and there could still be
non-unital inner and non-inner examples or we could be in characteristic 2.

In this paper we focus on this latter possibility for invariant metrics on
inner calculi by taking $k={\mathbb{F}}_2$ the field of two elements ${0,1}$. 
We take $\lambda=1$ and trivial $*$-structure as these do not make much
sense when there are only two elements. In the Abelian case condition (\ref
{metric}) becomes 
\begin{equation*}
B(v\circ w,z)=B(v,w\circ z)
\end{equation*}
and we can keep this also in the Heisenberg case of $B$ has its values in
the constants, which means the coefficients of the metric (given by $B^{-1}$
) are constants (an 'invariant metric'). In this case the condition on $B$
means that the data is precisely that of a commutative Frobenius algebra
over ${\mathbb{F}}_2$. So for each of these we obtain a differential
calculus and metric on the symmetric algebra on the vector space of the
algebra.

We close our generalities with a few general classes of examples:

(i) Let $X$ be a finite set of order $n$ and $V=k(X)$ the algebra of
functions on $X$ with pointwise product $\circ$. We let $x^\mu=\delta_\mu$ the delta
function at point $\mu\in X$ and since $x^\mu\circ x^\nu=\delta_{\mu\nu}x^\mu$ we have 
\begin{equation*}
[\mathrm{d} x^\mu,x^\nu]=\delta_{\mu\nu}\mathrm{d} x^\mu
\end{equation*}
with inner element $\theta=\sum_\mu \mathrm{d} x^\mu$. A for quantum metric $
g=\sum g_{\mu\nu}\mathrm{d} x^\mu\otimes\mathrm{d} x^\nu$ we require 
\begin{equation*}
[g,x^\rho]:=\sum g_{\mu\rho}\mathrm{d} x^\mu\otimes \mathrm{d} x^\rho+ \sum g_{\rho\mu}\mathrm{
d} x^\rho\otimes\mathrm{d} x^\mu=0
\end{equation*}
for all $\rho$ which implies $g_{\rho\rho}=0$ unless we are in characteristic 2 and $
g_{\mu\rho}=0$ for all $\mu\ne \rho$. So there is no metric unless we work in
characteristic 2 but in this case 
\begin{equation*}
g=\sum_\mu \mathrm{d} x^\mu\otimes \mathrm{d} x^\mu
\end{equation*}
is the unique quantum metric (the Euclidean metric) over $\F_2$.

(ii) $V=k{\mathbb{\ Z}}_n=k[x]/{\langle}x^n-e\>$ with basis $x^\mu$ the
different powers of $x$ with respect to $\circ$. We have commutation
relations 
\begin{equation*}
[\mathrm{d} x^\mu,x^\nu]=\mathrm{d} x^{\mu+\nu}
\end{equation*}
with indices treated mod $n$. At least $n$ quantum metrics exist over ${\mathbb{F}}_2$, namely the $n$ metrics
\begin{equation*}
g=\sum_\mu \mathrm{d} x^\mu\otimes \mathrm{d} x^{m-\mu};\quad  m=0,1,\cdots n-1.
\end{equation*}
We check that $[g,x^\rho]:=\sum \mathrm{d} x^{\mu+\rho}\otimes\mathrm{d} x^{m-\mu}+
\mathrm{d} x^\mu\otimes\mathrm{d} x^{m-\mu+\rho}=0$ after a relabelling $\mu+\rho\to \mu$
in the first sum to get 2 copies. In addition the elements $c=\sum_\mu \extd x^\mu$ and hence $c\tens c$ are central
and adding the latter gives a complementary metric where all the coefficients are reversed $0\leftrightarrow 1$. We will see that for
$n=2$ this gives no new metrics and indeed find just the above two, and for $n=3$ complementary metrics are degenerate so again give no more nondegenerate metrics and we find just the above 3 (this will not be the case for $n=4$ where we obtain 8 metrics). 

(iii) With $n=p^{d}$ and working over ${\mathbb{F}}_{p}$ where $p$ is prime,
there is a natural algebra $V=A_{d}=k[x]/{\langle }x^{p^{d}}-x \>$ which plays
an important role in the theory of field extensions. We have $x^\mu$ the
powers under $\circ$ with $e=x^0$ and $\mu=0,\cdots,n-1$. We focus on $p=2$.

$A_1$ is 2-dimensional with  $e$ a unit and $x\circ x=x$. The calculus is $[
\mathrm{d} e,e]=\mathrm{d} e$, $[\mathrm{d} e,x]=[\mathrm{d} x,e]=[\mathrm{d}
x,x]=\mathrm{d} x$. This is case B among the algebras for $n=2$ in the next section
and we find there that there is exactly one quantum metric $g= \mathrm{d}
e\otimes\mathrm{d} e+\mathrm{d} e\otimes\mathrm{d} x+\mathrm{d} x\otimes
\mathrm{d} e$. In fact this is isomorphic to (i) for 2 points.  $A_2$ is 4-dimensional with $e$ a unit, $x^\mu\circ x^\nu=x^{\mu+\nu}$
if $\mu+\nu<4$ and reduced by $x^4=x$ otherwise. Its own NCG was studied in \cite{BasMa} but now we are not
studying its NCG but rather that of $k[x^0,x^1,x^2,x^3]$ as a 4-dimensional
noncommutative spacetime. We will find 3 metrics for this calculus in Section~\ref{secn4}.

Once we have found the calculus and the metric we could hope to find a
quantum torsion free metric compatible or `quantum Levi-Civita' bimodule connecton (QLC for
short). By 'bimodule connection' on $\Omega ^{1}$ we mean a left connection,
i.e. $\nabla :\Omega ^{1}\rightarrow \Omega ^{1}\otimes _{A}\Omega ^{1}$
such that $\nabla (a\omega )=a(\nabla \omega )+da\otimes \omega $ for all $
a\in A,\omega \in \Omega ^{1}$ and in addition for some bimodule map $\sigma
:\Omega ^{1}\otimes _{A}\Omega ^{1}\rightarrow \Omega ^{1}\otimes _{A}\Omega
^{1}$:$\nabla (\omega a)=(\nabla \omega )a+\sigma (\omega \otimes da)$.

In \cite{Ma:gra} it is shown that in the inner case (with $\theta$) the
construction of a bimodule connection is equivalent to the construction of
bimodule maps $\sigma $ and $\alpha :\Omega ^{1}\rightarrow \Omega
^{1}\otimes _{A}\Omega ^{1}$. Then 
\begin{equation}  \label{ibc}
\nabla \omega =\theta \otimes \omega -\sigma \left( \omega \otimes \theta
\right) +\alpha \omega.
\end{equation}
Such a bimodule connection is metric compatible if
\begin{equation}  \label{imcbc}
\theta \otimes g+\left( \alpha \otimes \id\right) g+\sigma _{12}\left(
\id\otimes \left( \alpha -\sigma _{\theta }\right) \right) g=0
\end{equation}
where $\sigma _{\theta }=\sigma( (\ )\otimes\theta) $. This condition
results in quadratic relations for the coefficients of $\sigma$.

Finally, the curvature and torsion of a connection are
\begin{equation}\label{curv}
R_{\nabla } =\left( \extd\otimes \id-(\wedge\otimes \id) \left( \id\otimes \nabla
\right) \right) \nabla \quad : \Omega ^{1}\rightarrow \Omega ^{2}\otimes_A
\Omega^1
\end{equation}
\begin{equation*}
T_{\nabla }=\wedge \nabla -\extd \quad: \Omega ^{1}\rightarrow \Omega ^{2}
\end{equation*}
where $\wedge:\Omega ^{1}\otimes_A\Omega^1\rightarrow \Omega ^{2}$ is the
exterior product.
In our setting $\Omega$ over $A=\mathbb{F}_2[x^\mu]$ is generated by $\extd x^\mu$ with the relations: $\extd (x^\mu)^2=0$ and $\extd x^\mu \wedge \extd x^\nu =\extd x^\nu \wedge \extd x^\mu$.

In the inner case the construction of a torsion free bimodule connection is
equivalent \cite{Ma:gra} to the bimodule maps $\sigma $ and $\alpha $
satisfying
\begin{equation}  \label{itf}
\wedge \sigma =-\wedge \quad ,\quad \wedge \alpha =0.
\end{equation}

In order to solve these equations by computer we write out all the
conditions in terms of structure tensors starting with the pre-Lie algebra
in the form 
\begin{equation}
x^{\mu }\circ x^{\nu }=V_{\ \ \ \rho }^{\mu \nu }x^{\rho },\quad V_{\ \ \ \rho }^{\mu
\nu }\in {\mathbb{F}}_{2}
\end{equation}

For our polynomial or Heisenberg cases we need symmetry of the product so 
\begin{equation}
V_{\quad\rho }^{\mu \nu }-V_{{ \ \ \ }\rho }^{\nu \mu }=0
\label{cond1}
\end{equation}
and from (\ref{prelie2}) we need 
\begin{equation}
V_{{ \ \ }\lambda }^{\rho \nu }V_{{ \ \ \ }\gamma }^{\lambda \mu
}=V_{{ \ \ \ }\lambda }^{\rho \mu }V_{{ \ \ \ }\gamma }^{\lambda
\nu }
\end{equation}
which given commutativity is associativity of the product $\circ $ in this
case. For an inner calculus we have additionally 
\begin{equation*}
\theta \cdot V=\id,\quad \theta _{\mu }V_{{ \ \ \ }\rho }^{\mu \nu
}=\delta _{\rho }^{\nu }
\end{equation*}
for some vector $\theta =\theta _{\mu }\mathrm{d}x^{\mu }$, which means that 
$\theta $ is the identity for $\circ $.

The differential calculus induced by the pre-Lie algebra structure has the
commutation relations: 
\begin{equation}
\left[ \extd x^{\rho },x^{\nu }\right] =V_{{ \ \ \ }\mu }^{\rho \nu }\extd x^{\mu
},\quad V_{{ \ \ \ }\mu }^{\rho \nu }\in {\mathbb{F}}_{2}  \label{diffK}
\end{equation}

The conditions for quantum metric $g=g_{\mu \nu }dx^{\mu }\otimes dx^{\nu }$ 
$\in \Omega ^{1}\otimes _{A}\Omega ^{1}$ (which is quantum symmetric in the
sense $\wedge (g)=0$ and invertible in the sense of existence of a bimodule
map$(,):\Omega ^{1}\otimes _{A}\Omega ^{1}\rightarrow A$ such that $(\omega
,g_{1})g_{2}=\omega =g_{1}(g_{2},\omega )$ for all $\omega \in \Omega ^{1}$,
where $g=g_{1}\otimes g_{2}$ with sums of such terms understood) are 
\begin{equation}
g_{\lambda \nu }V_{{ \ \ \ }\mu }^{\lambda \rho }+g_{\mu \gamma }V_{
{ \ \ \ }\nu }^{\gamma \rho }=0
\end{equation}
where $g_{\mu \nu }=g_{\nu \mu }.$

For a bimodule connection (for inner calculi) we require bimodule maps $
\sigma $ and $\alpha :\Omega ^{1}\rightarrow \Omega ^{1}\otimes _{A}\Omega
^{1}$ as above. For alpha map on the basis 1-forms, taking $\alpha \left(
\extd x^{\mu }\right) =\alpha _{{ \ \ }\nu \rho }^{\mu }\extd x^{\nu }\otimes
\extd x^{\rho }$ and $\alpha _{{ \ \ }\nu \rho }^{\mu }=\alpha _{{ \ }
\rho \nu }^{\mu }$, we require

\begin{equation}\label{alphaK}
\alpha _{{ \ }\gamma \sigma }^{\rho }V_{{ \ \ }\lambda }^{\gamma
\nu }+\alpha _{{ \ }\lambda \gamma }^{\rho }V_{{ \ \ \ }\sigma
}^{\gamma \nu }=V_{{ \ \ \ }\mu }^{\rho \nu }\alpha _{{ \ \ }
\lambda \sigma }^{\mu }
\end{equation}

These conditions come from the compatibility of $\alpha $ with the
differential calculus, i.e. from equality $\alpha \left( \left[ \extd x^{\rho },x^{\nu }\right] \right) =V_{
{ \ \ \ }\mu }^{\rho \nu }\alpha \left( \extd x^{\mu }\right) $ calculating
the left hand side $\alpha \left( \left[ \extd x^{\rho },x^{\nu }\right] \right) =
\left[ \alpha \left( \extd x^{\rho }\right) ,x^{\nu }\right] =\alpha _{{ \ \ 
}\gamma \sigma }^{\rho }\left[ \extd x^{\gamma }\otimes \extd x^{\sigma },x^{\nu }
\right] =\alpha _{{ \ \ }\gamma \sigma }^{\rho }V_{{ \ \ \ }
\lambda }^{\gamma \nu }\extd x^{\lambda }\otimes \extd x^{\sigma }+$\\
$+\alpha _{{ \ \ 
}\lambda \gamma }^{\rho }V_{{ \ \ \ \ }\sigma }^{\gamma \nu
}\extd x^{\lambda }\otimes \extd x^{\sigma }$ and from the right hand side $V_{{ \ \ \ }\mu }^{\rho \nu }\alpha \left(
\extd x^{\mu }\right) $ $=$ $V_{{ \ \ \ }\mu }^{\rho \nu }\alpha _{{ \
\ }\lambda \sigma }^{\mu }\extd x^{\lambda }\otimes \extd x^{\sigma }$ gives the above
relation (\ref{alphaK}).

Similarly for the sigma map, we assume $\sigma \left( \extd x^{\mu }\otimes
\extd x^{\nu }\right) =\sigma _{_{{ \ \ \ \ }}\rho \lambda }^{\mu \nu
}\extd x^{\rho }\otimes \extd x^{\lambda }$ and require compatibility: $\sigma \left( 
\left[ \extd x^{\mu }\otimes \extd x^{\nu },x^{\gamma }\right] \right) =\left[ \sigma
\left( \extd x^{\mu }\otimes \extd x^{\nu }\right) ,x^{\gamma }\right] $. From the left hand side we obtain $\sigma \left( \left[ \extd x^{\mu }\otimes \extd x^{\nu
},x^{\gamma }\right] \right) =V^{\mu \gamma }{}_{{ \ \ }\alpha }\sigma
\left( \extd x^{\alpha }\otimes \extd x^{\nu }\right) +V_{{ \ \ \ }\beta }^{\nu
\gamma }\sigma \left( \extd x^{\mu }\otimes \extd x^{\beta }\right) =\left( V^{\mu
\gamma }{}_{{ \ \ }\alpha }\sigma _{{ \ \ \ }\omega \sigma
}^{\alpha \nu }+V_{{ \ \ \ }\beta }^{\nu \gamma }\sigma _{{ \ \ \ }
\omega \sigma }^{\mu \beta }\right) \extd x^{\omega }\otimes \extd x^{\sigma }$ and from
the right hand side we obtain\\ $\left[ \sigma \left( \extd x^{\mu }\otimes \extd x^{\nu }\right) ,x^{\gamma }
\right] =\left[ \sigma _{{ \ \ \ }\lambda \rho }^{\mu \nu }\extd x^{\lambda
}\otimes \extd x^{\rho },x^{\gamma }\right] =\left( \sigma _{{ \ \ \ }
\lambda \sigma }^{\mu \nu }V_{{ \ \ \ \ }\omega }^{\lambda \gamma
}\extd x^{\omega }\otimes \extd x^{\sigma }+\sigma _{{ \ \ \ }\omega \rho }^{\mu
\nu }V_{{ \ \ \ }\sigma }^{\rho \gamma }\extd x^{\omega }\otimes \extd x^{\sigma
}\right) $ which results in the condition 
\begin{equation}\label{sK}
V^{\mu \gamma }{}_{\alpha }\sigma _{{ \ \ \ \ }\omega \sigma }^{\alpha
\nu }+V_{{ \ \ \ }\beta }^{\nu \gamma }\sigma _{{ \ \ \ }\omega
\sigma }^{\mu \beta }=\sigma _{{ \ \ \ }\lambda \sigma }^{\mu \nu }V_{
{ \ \ \ \ }\omega }^{\lambda \gamma }+\sigma _{{ \ \ \ }\omega
\rho }^{\mu \nu }V_{{ \ \ \ }\sigma }^{\rho \gamma }
\end{equation}

Additionally we are interested in $\sigma $ invertible as $2n\times 2n$
matrix ${\mathbb{F}}_{2}^{2n}\rightarrow {\mathbb{F}}_{2}^{2n}$, i.e. 

\begin{equation*}
\sigma =\left( 
\begin{array}{cccc}
\sigma _{{ \ \ \ }11}^{11} & \sigma _{{ \ \ \ }12}^{11} & ... & 
\sigma _{{ \ \ \ }nn}^{11} \\ 
\sigma _{{ \ \ \ }11}^{12} &  &  &  \\ 
... &  &...  &  \\ 
\sigma _{{ \ \ \ }11}^{nn} &  &  & \sigma _{{ \ \ \ }nn}^{nn}
\end{array}
\right)\quad \mbox{with}\quad\det \left( \sigma \right) \neq 0.
\end{equation*}

Metric compatibility (\ref{imcbc}) (for $\alpha =0$, which turns out to be the
case in our considerations) with $\theta =\extd x^{\alpha }$  leads to the
non-linear conditions for the sigma coefficients $\sigma _{_{ \ \ \ }
\rho \lambda }^{\mu \nu }$,

\begin{equation}\label{ss}
g_{\gamma \rho }\delta _{\alpha \beta }-g_{\mu \nu }\sigma _{{ \ \ \ \ }
\lambda \rho }^{\nu \alpha }\sigma _{{ \ \ \ \ }\beta \gamma }^{\mu
\lambda }=0. 
\end{equation}

\section{Classification and their quantum geometries for $n=2$}\label{secn2}

For $n=1$ there are up to isomorphism only two algebras of dimension 1
namely $x\circ x=0$ which is nonunital and gives the classical calculus $[ 
\mathrm{d} x,x]=0$ on $k[x]$ and $e\circ e=e$ which gives the finite
difference calculus $[\mathrm{d} e,e]=\lambda \mathrm{d} e$ for deformation
parameter $\lambda$ (which we can take to be 1) in agreement with \cite
{Ma:fie} for 1-dimensional calculi. The only candidate for a quantum metric is $g=
\mathrm{d} x\otimes\mathrm{d} x$ which is central as the calculus is is
commutative and similarly $g=\mathrm{d} e\otimes\mathrm{d} e$ which is only
central over ${\mathbb{F}}_2$.

Our goal in this section is to give a full classification of the $n=2$
unital (hence inner) case. Clearly 2-dimensional commutative unital algebras
have the form $e,x$ as basis and $e\circ e=e,e\circ x=x=x\circ e$ and $
x\circ x=\alpha e+\beta x$ for constants $\alpha ,\beta $ for free
parameters $\alpha ,\beta $. Over ${\ \mathbb{F}}_{2}$ this means four
possibilities 
\begin{equation*}
x\circ x=0,\quad x\circ x=e,\quad x\circ x=x,\quad x\circ x=e+x
\end{equation*}
with the first two isomorphic by $x\mapsto x+e$. Thus there are three
inequivalent unital algebras hence three calculi. We use a computer to check
this explicitly (as a warm up to the next section), then find the metrics in
each case and the metric compatible quantum Levi-Civita connections. We summarise our results in Table 1.

\begin{figure}
\small{
\begin{tabular}{|l|l|l|l|l|}
\hline
& $
\begin{array}{c}
\text{Relations} \\ 
\left[ \extd e,e\right] =\extd e \\ 
\left[ \extd e,x\right] =\extd x=\left[ \extd x,e\right] 
\end{array}
$ 
 & {\tiny {Orbit order}} & {\tiny {Isotropy Group}} & Central
metrics and nonzero QLCs \\ \hline
A & $
\begin{array}{c}
x\circ x=0 \\ 
\left[ \mathrm{d}x,x\right] =0
\end{array}
$ & $2\times 3$ & $\{\mathfrak{1}\}$ & $
\begin{array}{c}
g_{A.I}=\mathrm{d}e\otimes \mathrm{d}x+\mathrm{d}x\otimes \mathrm{d}e \\ 
\nabla \mathrm{d}e=\mathrm{d}x\otimes \mathrm{d}x,\quad \nabla \mathrm{d}x=0
\\ 
\nabla \mathrm{d}e= \mathrm{d}x\otimes \mathrm{d}e+\mathrm{d}e\otimes 
\mathrm{d}x ,\quad \nabla \mathrm{d}x=\mathrm{d}x\otimes \mathrm{d}x
\\ 
\nabla \mathrm{d}e= \mathrm{d}x\otimes \mathrm{d}e+\mathrm{d}e\otimes 
\mathrm{d}x+\mathrm{d}x\otimes \mathrm{d}x ,\quad \nabla \mathrm{d}x=
\mathrm{d}x\otimes \mathrm{d}x \\ 
\nabla \mathrm{d}e=0,\quad \nabla \mathrm{d}x=\mathrm{d}e\otimes \mathrm{d}e
\\ 
\nabla \mathrm{d}e= \mathrm{d}e\otimes \mathrm{d}x+\mathrm{d}x\otimes 
\mathrm{d}e ,\quad \nabla \mathrm{d}x= \mathrm{d}e\otimes 
\mathrm{d}e+\mathrm{d}x\otimes \mathrm{d}x  \\ \hline
g_{A.II}=\mathrm{d}e\otimes \mathrm{d}x+\mathrm{d}x\otimes \mathrm{d}e+
\mathrm{d}x\otimes \mathrm{d}x \\ 
\nabla \mathrm{d}e=\mathrm{d}x\otimes \mathrm{d}x,\quad \nabla \mathrm{d}x=0
\\ 
\nabla \mathrm{d}e= \mathrm{d}e\otimes \mathrm{d}x+\mathrm{d}x\otimes 
\mathrm{d}e ,\quad \nabla \mathrm{d}x=\mathrm{d}x\otimes \mathrm{d}x
\\ 
\nabla \mathrm{d}e=\mathrm{d}e\otimes \mathrm{d}x+\mathrm{d}x\otimes 
\mathrm{d}e+\mathrm{d}x\otimes \mathrm{d}x ,\quad \nabla \mathrm{d}x=
\mathrm{d}x\otimes \mathrm{d}x
\end{array}
$ \\ \hline
B & $
\begin{array}{c}
x\circ x=x~ \\ 
\left[ \mathrm{d}x,x\right] =\mathrm{d}x
\end{array}
$ & $1\times 3$ & $\{\mathfrak{1},u\}=\mathbb{Z}_{2}$ & $
\begin{array}{c}
g_{B}=\mathrm{d}e\otimes \mathrm{d}e+\mathrm{d}e\otimes \mathrm{d}x+\mathrm{d
}x\otimes \mathrm{d}e \\ 
\nabla \mathrm{d}e=0,\quad \nabla \mathrm{d}x=\mathrm{d}e\otimes \mathrm{d}e
\end{array}
$ \\ \hline
C & $
\begin{array}{c}
x\circ x=e+x~ \\ 
\left[ \mathrm{d}x,x\right] =\mathrm{d}e+\mathrm{d}x
\end{array}
$ & $1\times 3$ & $\{\mathfrak{1},v\}=\mathbb{Z}_{2}$ & $
\begin{array}{c}
g_{C.I}=\mathrm{d}e\otimes \mathrm{d}x+\mathrm{d}x\otimes \mathrm{d}e+
\mathrm{d}x\otimes \mathrm{d}x \\ 
\nabla \mathrm{d}e= \mathrm{d}e\otimes \mathrm{d}x+\mathrm{d}x\otimes 
\mathrm{d}e+\mathrm{d}x\otimes \mathrm{d}x ,\quad \nabla \mathrm{d}x=
\mathrm{d}x\otimes \mathrm{d}x \\ \hline
g_{C.II}=\mathrm{d}e\otimes \mathrm{d}e+\mathrm{d}e\otimes \mathrm{d}x+
\mathrm{d}x\otimes \mathrm{d}e \\ 
\nabla \mathrm{d}e=0,\quad \nabla \mathrm{d}x=\mathrm{d}e\otimes \mathrm{d}e
\\ \hline
g_{C.III}=\mathrm{d}e\otimes \mathrm{d}e+\mathrm{d}x\otimes \mathrm{d}x \\ 
\nabla \mathrm{d}e=\mathrm{d}e\otimes \mathrm{d}e+\mathrm{d}x\otimes \mathrm{
d}x,\quad \nabla \mathrm{d}x=\mathrm{d}e\otimes \mathrm{d}x+\mathrm{d}
x\otimes \mathrm{d}e
\end{array}
$ \\ \hline
\end{tabular}
\newline
Table 1. All possible inner noncommutative geometries on $\F_2[e,x]$. Note that $g_{A.II}=g_{C.I}$ and $g_{B}=g_{C.II}$.}\normalsize\end{figure}

We now explain how these results were obtained. {\em A priori} the noncommutative geometry of interest is that of ${\mathbb{F}}_{2}[x^{1},x^{2}]$, defined as the universal enveloping
algebra of Abelian Lie algebra generated by basis elements $x^{1},x^{2}$,
the commutative algebra product (\ref{cond1}) $x^i\circ x^j=V^{ij}_{\quad l}x^l$ 
induces the differential calculus (\ref{diffK}). Notice that in 2 dimensions
with variables $x^{1}$ and $x^{2}$ we have three possibilities of inner
calculi with $\theta $ as the differential of an element of the pre-Lie
algebra and they are $\ \theta =\mathrm{d}x^{1}$ or $\mathrm{d}x^{2}$ or $
\mathrm{d}x^{1}+\mathrm{d}x^{2}$.

We find all the possible solutions of the commutative pre-Lie algebra
structure in 2 dimensions which induces inner differential calculus, i.e.
satisfying (\ref{prelie1}) and (\ref{prelie2}). Let $
S=\{s_{1},...,s_{12}\}$ be the set of all solutions and they can be grouped
as follows:

$\bullet$ 4 cases of inner calculus with $\theta =\mathrm{d}x^{1}$,
all have the following commutation relations $\left[ \mathrm{d}x^{1},x^{1}
\right] =\mathrm{d}x^{1},\left[ \mathrm{d}x^{1},x^{2}\right] =\mathrm{d}
x^{2}=\left[ \mathrm{d}x^{2},x^{1}\right] $ and the remaining commutators we
order as follows: $s_{1}:\left[ \mathrm{d}x^{2},x^{2}\right] =0$ , $s_{2}:
\left[ \mathrm{d}x^{2},x^{2}\right] =\mathrm{d}x^{2}$ , $s_{3}:\left[ 
\mathrm{d}x^{2},x^{2}\right] =\mathrm{d}x^{1}+\mathrm{d}x^{2}$ , $s_{4}:
\left[ \mathrm{d}x^{2},x^{2}\right] =\mathrm{d}x^{1}$;

$\bullet$ 4 cases of inner calculus with $\theta =\mathrm{d}x^{2}$
with $\left[ \mathrm{d}x^{2},x^{i}\right] =\mathrm{d}x^{i}=\left[ \mathrm{d}
x^{i},x^{2}\right] $ and $s_{5}:\left[ \mathrm{d}x^{1},x^{1}\right] =0$ , $
s_{6}:\left[ \mathrm{d}x^{1},x^{1}\right] =\mathrm{d}x^{2}$ , $s_{7}:\left[ 
\mathrm{d}x^{1},x^{1}\right] =\mathrm{d}x^{1}$ , $s_{8}:\left[ \mathrm{d}
x^{1},x^{1}\right] =\mathrm{d}x^{1}+\mathrm{d}x^{2}$;

$\bullet$ 4 cases of inner calculus with $\theta =\mathrm{d}x^{1}+
\mathrm{d}x^{2}$, such that $\left[ \mathrm{d}x^{1}+\mathrm{d}x^{2},x^{i}
\right] =\mathrm{d}x^{i}$ :

$s_{9}:\left[ \mathrm{d}x^{1},x^{1}\right] =0\quad ,\quad \left[ \mathrm{d}
x^{1},x^{2}\right] =\mathrm{d}x^{1}=\left[ \mathrm{d}x^{2},x^{1}\right]
\quad ,\quad \left[ \mathrm{d}x^{2},x^{2}\right] =\mathrm{d}x^{1}+\mathrm{d}
x^{2},$

$s_{10}:\left[ \mathrm{d}x^{1},x^{1}\right] =\mathrm{d}x^{2}\quad ,\quad 
\left[ \mathrm{d}x^{1},x^{2}\right] =\mathrm{d}x^{1}+\mathrm{d}x^{2}=\left[ 
\mathrm{d}x^{2},x^{1}\right] \quad ,\quad \left[ \mathrm{d}x^{2},x^{2}\right]
=\mathrm{d}x^{1},$

$s_{11}:\left[ \mathrm{d}x^{1},x^{1}\right] =\mathrm{d}x^{1}\quad ,\quad 
\left[ \mathrm{d}x^{1},x^{2}\right] =0=\left[ \mathrm{d}x^{2},x^{1}\right]
\quad ,\quad \left[ \mathrm{d}x^{2},x^{2}\right] =\mathrm{d}x^{2}$ ,

$s_{12}:\left[ \mathrm{d}x^{1},x^{1}\right] =\mathrm{d}x^{1}+\mathrm{d}
x^{2}\quad ,\quad \left[ \mathrm{d}x^{1},x^{2}\right] =\mathrm{d}x^{2}=\left[
\mathrm{d}x^{2},x^{1}\right] \quad ,\quad \left[ \mathrm{d}x^{2},x^{2}\right]
=0.$

One can show that due to the action of the group of isomorphisms $G$ on the
set $S$  of these solutions we get only three inequivalent families,
corresponding to the orbits of the action of the group.

The group of isomorphisms in 2 dimensions over $\mathbb{F}_{2}$ is $
G=SL\left( 2,2\right) =PSL\left( 2,2\right) =S_{3}$ (of order 6) with the
elements
\[ \mathfrak{1}=\left( 
\begin{array}{cc}
1 & 0 \\ 
0 & 1
\end{array}
\right) ,u=\left( 
\begin{array}{cc}
1 & 0 \\ 
1 & 1
\end{array}
\right) ,w=\left( 
\begin{array}{cc}
0 & 1 \\ 
1 & 0
\end{array}
\right) ,vu=\left( 
\begin{array}{cc}
0 & 1 \\ 
1 & 1
\end{array}
\right) ,uv=\left( 
\begin{array}{cc}
1 & 1 \\ 
1 & 0
\end{array}
\right) ,v=\left( 
\begin{array}{cc}
1 & 1 \\ 
0 & 1
\end{array}
\right) .\]
 Its action on the set of solutions $S$ results in the change of
variables, e.g. the action of the element $u$ corresponds to the 
change of variables 
\begin{equation}
\left( 
\begin{array}{c}
y^{1} \\ 
y^{2}
\end{array}
\right) =\left( 
\begin{array}{cc}
1 & 0 \\ 
1 & 1
\end{array}
\right) \left( 
\begin{array}{c}
x^{1} \\ 
x^{2}
\end{array}
\right) .  \label{change1}
\end{equation}

Note that already the set of the first $4$ solutions (for inner calculi with 
$\mathrm{d}x^{1}$,i.e. $S_{1}=\{s_{1},s_{2},s_{3},s_{4}\}$ splits into the
three orbits under the action of the group $G$. It is enough to consider the
element $u\in G$ and the change of variables (\ref{change1}) to obtain that $
s_{1}\simeq s_{4}$.

Recall that if $G$ acts on a set $S$ the orbits of this action are the sets 
\begin{equation*}
O_{s}=\{s^{\prime }\in S\quad |\quad g\cdot s=s^{\prime }\ \mbox{ for }\ g\in
G\}.
\end{equation*}

 We obtain the following:

For the calculus A the orbit consist of the elements: $O_{s_{1}}=
\{s_{1},s_{4,}s_{5},s_{6},s_{9},s_{12}\},|O_{s_{1}}|=6$ and the isotropy
group of element $Hs_{1}=\{\mathfrak{1}\}$.

For the calculus B: $O_{s_{2}}=\{s_{2},s_{7,}s_{11}\},\left\vert
O_{s_{2}}\right\vert =2$ and the isotropy group of the element $Hs_{2}=\{
\mathfrak{1},u\}$.

For the calculus C: $O_{s_{3}}=\{s_{3},s_{8,}s_{10}\},\left\vert
O_{s_{3}}\right\vert =2$ and the isotropy group of the element $Hs_{3}=\{
\mathfrak{1},v\}$.

As an example we present the explicit calculation for the orbit containing
element $s_{2}$: $e\cdot s_{2}=s_{2};u\cdot s_{2}=s_{2};$  $w\cdot s_{2}=s_{7};uv\cdot
s_{2}=s_{7};$  $vu\cdot s_{2}=s_{11};$  $v\cdot s_{2}=s_{11}$. Therefore the
corresponding orbit is 
\begin{equation*}
O_{s_{2}}=\{s_{2},s_{7,}s_{11}\}
\end{equation*}
and the isotropy groups of its elements are 
\begin{equation*}
H_{s_{2}}=\{e,u\},\quad H_{s_{7}}=\{w,uv\},\quad H_{s_{11}}=\{vu,v\}.
\end{equation*}

Other orbits are calculated analogously. These three orbits exhaust the elements of the whole set $S$ implying there
are only three non-isomorphic families of differential calculi as collected in the first column of 
Table 1, choosing $s_{1}$ as case A, $s_{2}$ as case B and $s_{3}$ as case C, with $x^{1}=e$ the
identity element for the $\circ$ product and $\theta =\mathrm{d}x^{1}=\mathrm{d}e$, and $x^{2}=x$.

\medskip
For each of these differential calculi A, B and C we next look for the central
metrics $g\in \Omega ^{1}\otimes _{A}\Omega ^{1}$ with $\wedge \left(
g\right) =0$ in the form
\begin{equation*}
g=g_{11}\mathrm{d}e\otimes \mathrm{d}e+g_{12}\left( \mathrm{d}e\otimes 
\mathrm{d}x+\mathrm{d}x\otimes \mathrm{d}e\right) +g_{22}\mathrm{d}x\otimes 
\mathrm{d}x
\end{equation*}
with constant coefficients, i.e. $g_{11},g_{12},g_{22}\in \mathbb{F}_{2}$.
Then we look for bimodule connections, which take the form (\ref{ibc})
including bimodule maps $\alpha $ and $\sigma $. Here $\nabla \mathrm{d} e=\nabla \mathrm{d} x=0$ and $\sigma ={\rm flip}$ on the
generators are always torsion free metric compatible bimodule connections but each of the calculi has an additional QLCs
which are collected along with the possible metrics in the last column of Table 1. 

We show some of the calculations behind the $A$ case explicitly, with similar arguments for the other cases. 
Thus, working with the $A$ calculus,  we first calculate: $\left[ g,e\right] =0$ and $\left[ g,x\right]
=g_{11}\left( \extd x\otimes \extd e+\extd e\otimes \extd x\right) =0\Rightarrow g_{11}=0$. The
possible (non degenerate) solutions for metric coefficients are 

i) $g_{12}=1,g_{22}=0$ resulting in $g_{A.I}=\extd e\otimes \extd x+\extd x\otimes \extd e$

ii) $g_{12}=1,g_{22}=1$ resulting in $g_{A.II}=\extd e\otimes \extd x+\extd x\otimes
\extd e+\extd x\otimes \extd x$.

Next, to find the bimodule connections, we look for the bimodule maps $
\alpha $ and $\sigma $, taking the former in the form $\alpha \left( \mathrm{d} x\right) =a\mathrm{d} e\otimes \mathrm{d}
e+b\left( \mathrm{d} e\otimes \mathrm{d} x+\mathrm{d} x\otimes \mathrm{d}
e\right) +c\mathrm{d} x\otimes \mathrm{d} x$ and we calculate that  $\alpha
\left( \left[ \mathrm{d} x,e\right] \right) =\left[ \alpha \left( \mathrm{d}
x\right) ,e\right] =  0$. On the other hand, for this calculus, $\alpha \left( \left[ \mathrm{d} x,e\right] \right)
=\alpha \left( \mathrm{d} x\right) =a\mathrm{d} e\otimes \mathrm{d}
e+b\left( \mathrm{d} e\otimes \mathrm{d} x+\mathrm{d} x\otimes \mathrm{d}
e\right) +c\mathrm{d} x\otimes \mathrm{d} x$. Therefore we have $a,b,c=0$.

Similarly for $\alpha \left( \mathrm{d} e\right) =a^{\prime }\mathrm{d}
e\otimes \mathrm{d} e+b^{\prime }\left( \mathrm{d} e\otimes \mathrm{d} x+
\mathrm{d} x\otimes \mathrm{d} e\right) +c^{\prime }\mathrm{d} x\otimes 
\mathrm{d} x$ we calculate that  $\alpha \left( \left[ \mathrm{d} e,e\right]
\right) =[\alpha \left( \mathrm{d} e\right) ,e]=0$, while on the other hand $\alpha
\left( \left[ \mathrm{d} e,e \right] \right) =\alpha \left( \mathrm{d}
e\right) =a^{\prime }\mathrm{d} e\otimes \mathrm{d} e+b^{\prime }\left( 
\mathrm{d} e\otimes \mathrm{d} x+\mathrm{d} x\otimes \mathrm{d} e\right)
+c^{\prime }\mathrm{d} x\otimes \mathrm{d} x$ implying $a^{\prime
},b^{\prime },c^{\prime }=0$. Hence there are no non-zero module maps $\alpha$.

For the sigma map we assume (\ref{sK}) and the metric
compatibility (with $\alpha =0$) (\ref{imcbc}). We solve the relations (\ref{sK}) and (\ref{ss}) over the field $\mathbb{F}_{2}$ by computer, which gives rise to the following torsion free metric
compatible `quantum Levi-Civita' bimodule connections. This gives us

i) For $g_{A.I}=\mathrm{d}e\otimes \mathrm{d}x+\mathrm{d}x\otimes 
\mathrm{d}e$ we have five solutions 

(A.I.1) $\nabla \mathrm{d}e=\mathrm{d}x\otimes \mathrm{d}x,\quad \nabla \mathrm{d}
x=0;$

(A.I.2) $\nabla \mathrm{d}e= \mathrm{d}x\otimes \mathrm{d}e+\mathrm{d}
e\otimes \mathrm{d}x ,\quad \nabla \mathrm{d}x=\mathrm{d}x\otimes 
\mathrm{d}x;$

(A.I.3) $\nabla \mathrm{d}e= \mathrm{d}x\otimes \mathrm{d}e+\mathrm{d}
e\otimes \mathrm{d}x+\mathrm{d}x\otimes \mathrm{d}x ,\quad \nabla 
\mathrm{d}x=\mathrm{d}x\otimes \mathrm{d}x;$

(A.I.4) $\nabla \mathrm{d}e=0,\quad \nabla \mathrm{d}x=\mathrm{d}e\otimes \mathrm{d}
e;$

(A.I.5) $\nabla \mathrm{d}e= \mathrm{d}e\otimes \mathrm{d}x+\mathrm{d}
x\otimes \mathrm{d}e ,\quad \nabla \mathrm{d}x= \mathrm{d}
e\otimes \mathrm{d}e+\mathrm{d}x\otimes \mathrm{d}x $.\newline

ii) For $g_{A.II}=\mathrm{d}e\otimes \mathrm{d}x+\mathrm{d}x\otimes \mathrm{d}
e+\mathrm{d}x\otimes \mathrm{d}x$ we have three solutions

(A.II.1) $\nabla \mathrm{d}e=\mathrm{d}x\otimes \mathrm{d}x,\quad \nabla \mathrm{d}
x=0;$

(A.II.2) $\nabla \mathrm{d}e= \mathrm{d}e\otimes \mathrm{d}x+\mathrm{d}
x\otimes \mathrm{d}e ,\quad \nabla \mathrm{d}x=\mathrm{d}x\otimes 
\mathrm{d}x;$

(A.II.3) $\nabla \mathrm{d}e= \mathrm{d}e\otimes \mathrm{d}x+\mathrm{d}
x\otimes \mathrm{d}e+\mathrm{d}x\otimes \mathrm{d}x ,\quad \nabla 
\mathrm{d}x=\mathrm{d}x\otimes \mathrm{d}x$.

We solve the $B,C$ cases similarly, all results being collected in
Table 1 above. We  have not listed the associated $\sigma$ as these are uniquely determined by  $\nabla$ and the commutation
relations.

As a check we see that calculus case A corresponds to $V\backsimeq \mathbb{F}_2\mathbb{Z}_2$ as described in the general analysis, see example (ii), in Section~\ref{secformalism}. The metrics $g_{A.II},g_{A.I}$ recover the two metrics there for $m=0,1$ after the change of variables to $x^1=e+x$, $x^0=e$ (where the superscript on the left is a label not an exponent). The complementary metrics in the general analysis duplicate these. The calculus case B corresponds to $V\backsimeq \mathbb{F}_2[\text{2 points}]$ in the general analysis, example (i) in Section~\ref{secformalism}), with the  metric $g_B$ agreeing with the Euclidean metric there on change of variables $x^0=e+x$, $x^1=x$. 

\begin{proposition}\label{flat2dim}
For $n=2$ all quantum Levi-Civita connections as listed in Table 1 are flat.
\end{proposition}
\proof
Explicitly using \eqref{curv} we demonstrate the calculation for
the first two bimodule connections compatible with the first metric $g_{A.I}$ in the family A. For (A.I.1) this is

 $R_{\nabla }\mathrm{d}e=\left( \mathrm{d}\otimes \id-\wedge \left( \id\otimes
\nabla \right) \right) \nabla \mathrm{d}e=-\wedge \left( \id\otimes \nabla
\right) \left( \mathrm{d}x\otimes \mathrm{d}x\right) =-\mathrm{d}x\wedge
\nabla \mathrm{d}x=0$,

$R_{\nabla }\mathrm{d}x=\left( \mathrm{d}\otimes \id-\wedge \left(
\id\otimes \nabla \right) \right) \nabla \mathrm{d}x=0$,

while for (A.I.2) the calculation is

 $R_{\nabla }\mathrm{d}e=-\wedge \left( \id\otimes \nabla \right) \left( 
\mathrm{d}x\otimes \mathrm{d}e+\mathrm{d}e\otimes \mathrm{d}x\right) =
\mathrm{d}x\wedge \nabla \mathrm{d}e+\mathrm{d}e\wedge \nabla \mathrm{d}x=
\mathrm{d}x\wedge (\mathrm{d}x\otimes \mathrm{d}e)+\mathrm{d}x\wedge (\mathrm{d}
e\otimes \mathrm{d}x)+\mathrm{d}e\wedge (\mathrm{d}x\otimes \mathrm{d}x)=0$,

$R_{\nabla }\mathrm{d}x=-\wedge \left( \id\otimes \nabla \right)
\left( \mathrm{d}x\otimes \mathrm{d}x\right) =-\mathrm{d}x\wedge \nabla 
\mathrm{d}x=\mathrm{d}x\wedge (\mathrm{d}x\otimes \mathrm{d}x)=0$.

Similarly for all the other bimodule connections
in Table 1 above. We refer only to bimodule quantum Levi-Civita connections with constant coefficients and invertible $\sigma$ as listed in the table. \endproof

\section{Classification for $n=3$}\label{secn3}

For the $n=3$ inner case over ${\mathbb{F}}_{2}$ we will find six
inequivalent unital algebras. In each case we take $e,x,y$ as basis and have $
e\circ e=e,e\circ x=x=x\circ e,e\circ y=y=y\circ e$, with the remaining
relations as:

A:\qquad $x\circ y=0=y\circ x,\quad x\circ x=0=y\circ y,$

B:\qquad $x\circ y=0=y\circ x,\quad x\circ x=x,\quad y\circ y=y,$

C:\qquad $x\circ y=0=y\circ x,\quad x\circ x=x,\quad y\circ y=0,$

D:\qquad $x\circ y=x+y=y\circ x,\quad x\circ x=y,\quad y\circ y=x,$

E:\qquad $x\circ y=0=y\circ x,\quad x\circ x=y,\quad y\circ y=0,$

F:\qquad $x\circ y=x+y=y\circ x,\quad x\circ x=e+x+y,\quad y\circ y=x.$

These six inequivalent commutative unital algebras imply six noncommutative
differential calculi as shown in Table 2. For each of them we show the number of quantum metrics and for each metric the number of torsion free cotorsion free (`quantum Levi-Civita') bimodule connections. The metrics are listed in detail in Table~4.\\

\begin{figure}\small{
\begin{tabular}{|l|l|l|l|l|l|l|}
\hline
& $
\begin{array}{c}
\text{Relations} \\ 
\left[ \extd e,e\right] =\extd e \\ 
\left[ \extd e,x\right] =\extd x=\left[ \extd x,e\right]  \\ 
\left[ \extd e,y\right] =\extd y=\left[ \extd y,e\right] 
\end{array}
$ & Orbit order & Isotropy Group & Metrics
& 
 
$\begin{array}{c}
\text{Nonzero}\\
 \text{QLC}
\end{array}$
& $\begin{array}{c}
 R_\nabla\neq 0
\end{array}$
\\ \hline
A & $
\begin{array}{c}
\left[ \extd x,y\right] =0=\left[ \extd y,x\right]  \\ 
\left[ \extd x,x\right] =0=\left[ \extd y,y\right] 
\end{array}
$ & $\left\vert O_{s_{1}}\right\vert =4\times 7$ & $
\begin{array}{c}
\{\mathfrak{1},\tilde{w},\tilde{u}\tilde{v},\tilde{v},\tilde{u},\tilde{v}\tilde{u}\}=S_3
\end{array}
$ & $0$ & 
- & -\\ \hline
B & $
\begin{array}{c}
\left[ \extd x,y\right] =0=\left[ \extd y,x\right]  \\ 
\left[ \extd x,x\right] =\extd x,\left[ \extd y,y\right] =\extd y
\end{array}
$ & $\left\vert O_{s_{34}}\right\vert =4\times 7$ & $
\begin{array}{c}
\{\mathfrak{1},\tilde{w},\tilde{\tilde{v}}\tilde{\tilde{u}},\tilde{\tilde{u}},\tilde{\tilde{v}},\tilde{\tilde{u}}\tilde{\tilde{v}}\}=S_3
\end{array}
$ & $1$ & 
$3$& $0$\\ \hline
C & $~
\begin{array}{c}
\left[ \extd x,y§\right] =0=\left[ \extd y,x\right]  \\ 
\left[ \extd x,x\right] =\extd x,\left[ \extd y,y\right] =0
\end{array}
$ & $\left\vert O_{s_{2}}\right\vert =24\times 7$ & $
\begin{array}{c}
\{\mathfrak{1}\}
\end{array}
$ & $2$ & {$
\begin{array}{c}
13\ \\ \text{each}\ g_{C} 
\end{array}$}
&
\tiny{$
\begin{array}{c}
2 \text{ for } g_{C.I}\\
3 \text{ for } g_{C.II}
\end{array}$}
\\ \hline
D & $
\begin{array}{c}
\left[ \extd x,y\right] =\extd x+\extd y=\left[ \extd y,x\right]  \\ 
\left[ \extd x,x\right] =\extd y,\left[ \extd y,y\right] =\extd x
\end{array}
$ & $\left\vert O_{s_{23}}\right\vert =12\times 7$ & $
\begin{array}{c}
\{\mathfrak{1},\tilde{w}\}=\Z_2
\end{array}
$ & $3$ & {$
\begin{array}{c}
3 \\ \text{each}\ g_{D} 
\end{array}$}
&
\tiny{$
\begin{array}{c}
1 \text{ for } g_{D.I}\\
0 \text{ for } g_{D.II}\\
1 \text{ for } g_{D.III}
\end{array}$}
\\ \hline
E & $~
\begin{array}{c}
\left[ \extd x,y\right] =0=\left[ \extd y,x\right]  \\ 
\left[ \extd x,x\right] =\extd y,\left[ \extd y,y\right] =0
\end{array}
$ & $\left\vert O_{s_{3}}\right\vert =12\times 7$ & $
\begin{array}{c}
\{\mathfrak{1},\tilde{v}\}=\Z_2
\end{array}
$ & $4$ & {$
\begin{array}{c}
13 \\ \text{each}\ g_{E} 
\end{array}$}&
\tiny{$
\begin{array}{c}
2 \text{ for } g_{E.I}\\
3 \text{ for } g_{E.II}\\
5 \text{ for } g_{E.III}\\
4 \text{ for } g_{E.IV}
\end{array}$}
\\ \hline
F & $~
\begin{array}{c}
\left[ \extd x,y\right] =\extd e+\extd x=\left[ \extd y,x\right]  \\ 
\left[ \extd x,x\right] =\extd e+\extd x+\extd y\\
\left[ \extd y,y\right] =\extd x
\end{array}
$ & $\left\vert O_{s_{20}}\right\vert =8\times 7$ & $
\begin{array}{c}
\{\mathfrak{1},\tilde{\tilde{v}}\tilde{\tilde{u}},\tilde{\tilde{u}}\tilde{\tilde{v}}\}=\Z_3
\end{array}
$ & $7$ & {$
\begin{array}{c}
3 \\ \text{each}\ g_{F} 
\end{array}$}&
\tiny{$
\begin{array}{c}
2 \text{ for each } g_{F}\\ 
\text{ except }\\
 0 \text{ for } g_{F.II}
\end{array}$}
\\ \hline
\end{tabular}
\newline
Table 2. All possible inner noncommutative geometries on $\F_2[e,x,y]$.}\medskip
\normalsize\end{figure}

We now outline how these results were obtained. For ${\mathbb{F}}_{2}[x^{1},x^{2},x^{3}]$ defined as universal enveloping
algebra of Abelian Lie algebra generated by the basis elements $
x^{1},x^{2},x^{3},$ the commutative product (\ref{cond1}) induces the differential calculus, as before. In 3 dimensions there is seven possibilities for element $\theta $ (as the
differential of an element of the pre-Lie algebra) for inner calculi, namely

$\theta =\mathrm{d}x^{1},\mathrm{d}x^{2},\mathrm{d}x^{3},\mathrm{d}x^{1}+
\mathrm{d}x^{2},\mathrm{d}x^{1}+\mathrm{d}x^{3}$, $\mathrm{d}x^{2}+\mathrm{d}x^{3}$ and $\mathrm{d}x^{1}+
\mathrm{d}x^{2}+\mathrm{d}x^{3}$.

Finding explicitly the solutions to (\ref{prelie1}), (\ref{prelie2}) gives $7\times 64$ cases, i.e. 64 solutions for each of the seven
possible $\theta $.  The 64 cases with inner calculus with $\theta=\mathrm{d}x^{1}$ are listed in Table~3, there are similarly 64 for each of the other 6 cases.

\begin{figure}
\[ \includegraphics[scale=0.75]{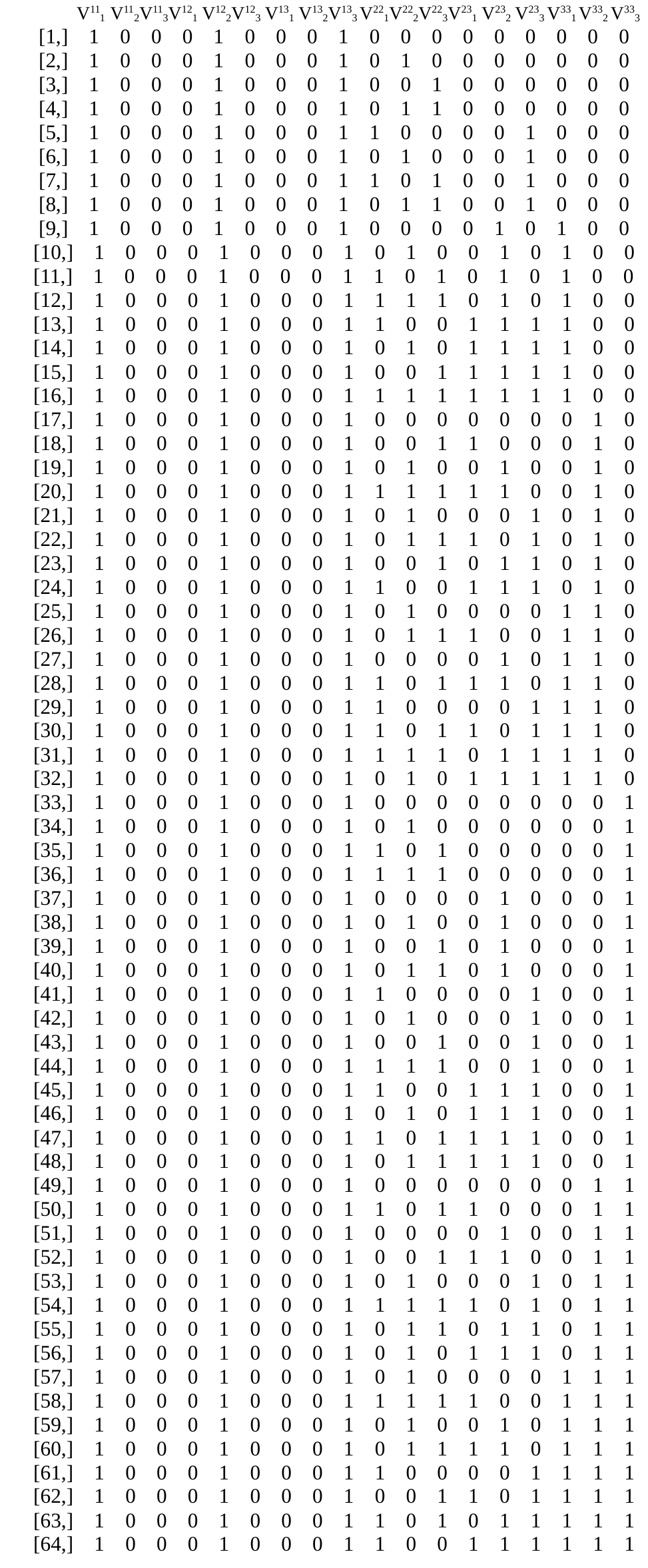}\]
Table~3. Structure coefficients for inner differential calculi with $\theta=\extd x^1$
\end{figure}

The set of the $7\times 64$ solutions, denoted by $S=
\{s_{1},...,s_{64},......s_{7\times 64}\}$, splits into the six orbits under
the action of a group of isomorphisms over ${\mathbb{F}}_{2}$: $G=GL\left( 3,
{\mathbb{F}}_{2}\right) =PSL\left( 2,7\right) $. The order of the group of
isomorphisms is $|G|=168$ and we write only some of its elements (only the
ones needed to list the isotropy groups explicitly in Table 2) :

\begin{eqnarray*}
\mathfrak{1} &=&\left( 
\begin{array}{ccc}
1 & 0 & 0 \\ 
0 & 1 & 0 \\ 
0 & 0 & 1
\end{array}
\right),\
 \tilde{u}=\left( 
\begin{array}{ccc}
1 & 0 & 0 \\ 
0 & 1 & 0 \\ 
0 & 1 & 1
\end{array}
\right),\ 
\tilde{w}=\left( 
\begin{array}{ccc}
1 & 0 & 0 \\ 
0 & 0 & 1 \\ 
0 & 1 & 0
\end{array}
\right) ,\ 
\tilde{u}\tilde{v}=\left( 
\begin{array}{ccc}
1 & 0 & 0 \\ 
0 & 1 & 1 \\ 
0 & 1 & 0
\end{array}
\right) ,\ 
\tilde{v}=\left( 
\begin{array}{ccc}
1 & 0 & 0 \\ 
0 & 1 & 1 \\ 
0 & 0 & 1
\end{array}
\right) , \\
 \tilde{v}\tilde{u}&=&\left( 
\begin{array}{ccc}
1 & 0 & 0 \\ 
0 & 0 & 1 \\ 
0 & 1 & 1
\end{array}
\right),\ 
\tilde{\tilde v}\tilde{\tilde u}=\left( 
\begin{array}{ccc}
1 & 0 & 0 \\ 
1 & 1 & 1 \\ 
0 & 1 & 0
\end{array}
\right),\ 
\tilde{\tilde u}=\left( 
\begin{array}{ccc}
1 & 0 & 0 \\ 
1 & 1 & 1 \\ 
0 & 0 & 1
\end{array}
\right),\ 
\tilde{\tilde v }=\left( 
\begin{array}{ccc}
1 & 0 & 0 \\ 
0 & 1 & 0 \\ 
1 & 1 & 1
\end{array}
\right),\ 
\tilde{\tilde u}\tilde{\tilde v}=\left( 
\begin{array}{ccc}
1 & 0 & 0 \\ 
0 & 0 & 1 \\ 
1 & 1 & 1
\end{array}
\right). 
\end{eqnarray*}

One can present the sketch of a proof, by considering only the solutions for
one of the possible $\theta$, for which we take $\theta =\mathrm{d}x^{1}
$. These are the first 64 solutions we denote by $
S_{1}=\{s_{1},...,s_{64}\}$ listed in each row of Table~3,  where coefficients $V_{
\text{ \ \ }\rho }^{\mu \nu }$ (as solutions of (\ref{prelie1}), (\ref{prelie2})  are collected. 
The other choices of $\theta$ have  similar structure by a $GL(3,\F_2)$ transformation mapping 
$\extd x^1$ to any other $\theta$. The space $S_1$ of such restricted solutions again splits into the six orbits under the action of
 $GL\left( 3,{\mathbb{F}}_{2}\right) $.

Thus,  $s_{i}$ corresponds to the case with coefficients $V_{\text{ \ \ }
\rho }^{\mu \nu }$ listed in the row $[i]$ in Table 3 and falls to the orbit denoted as $
O_{s_{i}}\cap O_{r_{i}}$ (which is equivalent to the one of the families A - F given in the Table 2 above). We write $O_{s_{i}}\cap O_{p_{i}}$to underline that we list
only part of the orbit (the one coming from the first 64 solutions with
inner calculus $\mathrm{d}x^{1}$ as $O_{p_{i}}$, when in fact the full size
of $O_{s_{i}}$ is 7 $\times $bigger, as indicated in Table 2).

The first 64 solutions split into the following orbits:

$O_{s_{1}}\cap O_{p_{1}}=\{s_{1},s_{5},s_{9},s_{13}\};$

$O_{s_{34}}\cap O_{p_{34}}=\{s_{34},s_{38},s_{42},s_{46}\};$

$O_{s_{2}}\cap
O_{p_{2}}=
\{s_{2},s_{4},s_{6},s_{8},s_{10},s_{12},s_{14},s_{16},s_{19},s_{21},s_{25},s_{32},s_{33},s_{35},s_{37},s_{39},s_{41},s_{43},s_{45},s_{47},s_{49},s_{51},s_{61},s_{64}\};
$

$O_{s_{23}}\cap O_{p_{23}}=\
\{s_{23},s_{18},s_{28},s_{30},s_{36},s_{40},s_{44},s_{48},s_{53},s_{56},s_{57},s_{59}\};
$

$O_{s_{3}}\cap
O_{p_{3}}=
\{s_{3},s_{7},s_{11},s_{15},s_{17},s_{24},s_{27},s_{29},s_{54},s_{55},s_{58},s_{60}\};
$

$O_{s_{20}}\cap
O_{p_{20}}=\{s_{20},s_{22},s_{26},s_{31},s_{50},s_{52},s_{62},s_{63}\}$.

Already for the first 64 solutions we get the six orbits giving six (A -
F) inequivalent noncommutative differential calculi.

For the isotropy groups one can calculate that, e.g. for $H_{s_{1}}$: $
\mathfrak{1} \cdot s_{1}=s_{1};\tilde{w}\cdot s_{1}=s_{1};\tilde{u}\tilde{v}\cdot
s_{1}=s_{1};\tilde{v}\cdot s_{1}=s_{1};\tilde{u}\cdot s_{1}=s_{1};\tilde{v}\tilde{u}\cdot
s_{1}=s_{1}$ or for $H_{s_{3}}$: $\mathfrak{1} \cdot s_{3},\tilde{v}\cdot s_{3}$ or
for $H_{s_{23}}$: $\mathfrak{1} \cdot s_{23}=s_{23}$ and $\tilde{w}\cdot s_{23}=s_{23}$. Similarly one can calculate the isotropy groups for the remaining elements.

Below we also show some examples of isomorphisms of certain solutions in $
S_{1}$ to the six presented above families.

One element of $G$ namely $\left( 
\begin{array}{ccc}
1 & 0 & 0 \\ 
0 & 0 & 1 \\ 
1 & 1 & 0
\end{array}
\right) $ gives the  change of variables 
\begin{eqnarray*}
y^{1} &=&x^{1};\quad y^{2}=x^{3};\quad y^{3}=x^{1}+x^{2} \\
\extd y^{1} &=&\extd x^{1};\quad \extd y^{2}=\extd x^{3};\quad \extd y^{3}=\extd x^{1}+\extd x^{2}
\end{eqnarray*}
under which we immediately see that 
\begin{equation*}
\left[ \extd y^{1},y^{1}\right] =\extd y^{1};\quad \left[ \extd y^{1},y^{2}\right] =\extd y^{2}=
\left[ \extd y^{2},y^{1}\right] ;\quad \left[ \extd y^{1},y^{3}\right] =\extd y^{3}=\left[
\extd y^{3},y^{1}\right] 
\end{equation*}
i.e. the result is still an inner differential calculus with $\theta
=\extd y^{1}=\extd x^{1}$.

One can check that the remaining commutators will always fall into one of
the six families above A - F:

$\bullet$ In case A, $s_{1}$ is isomorphic to $s_{9}$, with 
$\left[ \extd y^{2},y^{2}\right] =0;\left[ \extd y^{2},y^{3}\right] =\extd y^{2}=\left[
\extd y^{3},y^{2}\right] ;\left[ \extd y^{3},y^{3}\right] =\extd y^{1}$ i.e. with the
non-zero coefficients, $V^{23}_{\quad 2}=V^{32}_{\quad 2}=1$ and $V^{33}_{\quad 1}=1$, cf. row [9] in \
Table~3 above.

$\bullet$ In case B,  $s_{34}$ is isomorphic to $s_{38}$, with $
\left[ \extd y^{2},y^{2}\right] =\extd y^{2};\left[ \extd y^{2},y^{3}\right] =\extd y^{2}=\left[
\extd y^{3},y^{2}\right] ;\left[ \extd y^{3},y^{3}\right] =\extd y^{3}$.

$\bullet$ In case C, $s_{2}$ is isomorphic to $s_{37}$, with
$\left[ \extd y^{2},y^{2}\right] =0;\left[ \extd y^{2},y^{3}\right] =\extd y^{2}=\left[
\extd y^{3},y^{2}\right] ;\left[ \extd y^{3},y^{3}\right] =\extd y^{3}$.

$\bullet$ In case D, $s_{23}$ is isomorphic to $s_{30}$, with  $\left[ \extd y^{2},y^{2}\right] =\extd y^{1}+\extd y^{3};\left[ \extd y^{2},y^{3}\right]
=\extd y^{1}+\extd y^{3}=\left[ \extd y^{3},y^{2}\right] ;\left[ \extd y^{3},y^{3}\right]
=\extd y^{1}+\extd y^{2}$.

$\bullet$ In case E, $s_{3}$ is isomorphic to $s_{27}$, with
$\left[ \extd y^{2},y^{2}\right] =0;\left[ \extd y^{2},y^{3}\right] =\extd y^{2}=\left[
\extd y^{3},y^{2}\right] ;\left[ \extd y^{3},y^{3}\right] =\extd y^{1}+\extd y^{2}$. 

$\bullet$ In case F, $s_{20}$ is isomorphic to $s_{63}$, with 
$\left[ \extd y^{2},y^{2}\right] =\extd y^{1}+\extd y^{3};\left[ \extd y^{2},y^{3}\right]
=\extd y^{2}+\extd y^{3}=\left[ \extd y^{3},y^{2}\right] ;\left[ \extd y^{3},y^{3}\right]
=\extd y^{1}+\extd y^{2}+\extd y^{3}$.

\subsection{Central metrics}

For each of the differential calculi A - F  in the first column of Table 2 we next look for the quantum 
metrics $g\in \Omega ^{1}\otimes _{A}\Omega ^{1}$ with $\wedge \left(
g\right) =0$ in the form  
\begin{eqnarray}
g&=&g_{11}\mathrm{d}e\otimes \mathrm{d}e+g_{12}\left( \mathrm{d}e\otimes 
\mathrm{d}x+\mathrm{d}x\otimes \mathrm{d}e\right) 
+g_{13}\left( \mathrm{d}e\otimes 
\mathrm{d}y+\mathrm{d}y\otimes \mathrm{d}e\right)
+g_{22}\mathrm{d}x\otimes 
\mathrm{d}x\nonumber\\
&+&g_{23}\left( \mathrm{d}x\otimes 
\mathrm{d}y+\mathrm{d}y\otimes \mathrm{d}x\right) +g_{33}\mathrm{d}y\otimes 
\mathrm{d}y\nonumber
\end{eqnarray}
with constant coefficients, i.e. $g_{\mu\nu}\in \mathbb{F}_{2}$, similarly as outlined in the $n=2$ case. The results for the 
possible  quantum metrics corresponding to the above differential calculi are shown in Table~4. 

\begin{figure}\small {
\begin{tabular}{|l|l|}
\hline
& Central metrics \\ \hline
A & 0 \\ \hline
B & $g_{B}=\extd e\otimes \extd e+\extd e\otimes \extd x+\extd x\otimes \extd e+\extd e\otimes \extd y+\extd y\otimes
\extd e+\extd x\otimes \extd y+\extd y\otimes \extd x$ \\ \hline
C & $
\begin{array}{c}
g_{C.I}=\extd e\otimes \extd y+\extd y\otimes \extd e+\extd x\otimes \extd x+\extd x\otimes \extd y+\extd y\otimes \extd x \\ 
g_{C.II}=\extd e\otimes \extd y+\extd y\otimes \extd e+\extd x\otimes \extd x+\extd x\otimes \extd y+\extd y\otimes
\extd x+\extd y\otimes \extd y
\end{array}
$ \\ \hline
D & $
\begin{array}{c}
g_{D.I}=\extd e\otimes \extd e+\extd e\otimes \extd x+\extd x\otimes \extd e+\extd e\otimes \extd y+\extd y\otimes
\extd e+\extd x\otimes \extd x \\ 
g_{D.II}=\extd e\otimes \extd e+\extd e\otimes \extd x+\extd x\otimes \extd e+\extd e\otimes \extd y+\extd y\otimes
\extd e+\extd x\otimes \extd y+\extd y\otimes \extd x \\ 
g_{D.III}=\extd e\otimes \extd e+\extd e\otimes \extd x+\extd x\otimes \extd e+\extd e\otimes \extd y+\extd y\otimes
\extd e+\extd y\otimes \extd y
\end{array}
$ \\ \hline
E & $
\begin{array}{c}
g_{E.I}=\extd e\otimes \extd y+\extd y\otimes \extd e+\extd x\otimes \extd x \\ 
g_{E.II}=\extd e\otimes \extd y+\extd y\otimes \extd e+\extd x\otimes \extd x+\extd x\otimes \extd y+\extd y\otimes \extd x \\ 
g_{E.III}=\extd e\otimes \extd y+\extd y\otimes \extd e+\extd x\otimes \extd x+\extd y\otimes \extd y \\ 
g_{E.IV}=\extd e\otimes \extd y+\extd y\otimes \extd e+\extd x\otimes \extd x+\extd x\otimes \extd y+\extd y\otimes
\extd x+\extd y\otimes \extd y
\end{array}
$ \\ \hline
F & $
\begin{array}{c}
g_{F.I}=\extd e\otimes \extd y+\extd y\otimes \extd e+\extd x\otimes \extd x \\ 
g_{F.II}=\extd e\otimes \extd e+\extd e\otimes \extd x+\extd x\otimes \extd e+\extd e\otimes \extd y+\extd y\otimes
\extd e+\extd x\otimes \extd y+\extd y\otimes \extd x \\ 
g_{F.III}=\extd e\otimes \extd e+\extd e\otimes \extd x+\extd x\otimes \extd e+\extd x\otimes \extd x+\extd x\otimes
\extd y+\extd y\otimes \extd x \\ 
g_{F.IV}=\extd e\otimes \extd x+\extd x\otimes \extd e+\extd e\otimes \extd y+\extd y\otimes \extd e+\extd y\otimes \extd y \\ 
g_{F.V}=\extd e\otimes \extd x+\extd x\otimes \extd e+\extd x\otimes \extd x+\extd y\otimes \extd y \\ 
g_{F.VI}=\extd e\otimes \extd e+\extd x\otimes \extd y+\extd y\otimes \extd x+\extd y\otimes \extd y \\ 
g_{F.VII}=\extd e\otimes \extd e+\extd e\otimes \extd y+\extd y\otimes \extd e+\extd x\otimes
\extd x+\extd x\otimes \extd y+\extd y\otimes \extd x+\extd y\otimes \extd y
\end{array}
$ \\ \hline
\end {tabular}\\
Table~4. All possible quantum metrics for each of the possible calculi in Table~2. Note that $g_{E.II}=g_{C.I}$, $g_{E.IV}=g_{C.II}$, $g_{F.I}=g_{E.I}$ and $g_{F.II}=g_{D.II}=g_B$. 
} \normalsize
\end{figure}

As a check we see that the calculus case B corresponds to $V\backsimeq \mathbb{F}_2[\text{3 points}]$ with the metric $g_{B}$ agreeing with the Euclidean metric in the general analysis (example (i) from Section~\ref{secformalism} after the change of variables: $e=x^0+x^1+x^2$, $x^1=x$, $x^2=y$. The case D corresponds to $V\backsimeq \mathbb{F}_2\mathbb{Z}_3$ with the metrics $g_{D.II}$, $g_{D.III},g_{D.I}$  recovering the  $m=0,1,2$  metrics for example (ii) in Section~\ref{secformalism}  after the change of variables $x^0=e, x^1=e+x, x^2=e+y$  (where the superscript of $x$ is a label not a square; this element being $(x^1)^2$ in $\F_2\Z_3$). The complementary metrics are not listed in our table as they are degenerate. 

\medskip
Next we look for bimodule connections, which take the form (\ref{ibc})
including bimodule maps $\alpha $ and $\sigma $. As for $n=2$, careful analysis shows that there are no nonzero module maps $\alpha$. For the $\sigma$ map we assume metric compatibility (\ref{imcbc}) and impose the torsion free condition (\ref{itf}). The curvature is calculated from \eqref{curv}. $\nabla \mathrm{d} e=\nabla \mathrm{d} x=\nabla \mathrm{d} y=0$ and $\sigma ={\rm flip}$ on the
generators are always torsion free metric compatible bimodule connections but we also find additional non-zero QLCs, some of which have non-zero curvature $R_\nabla$ (their numbers are given in the Table 2). This methodology is the same as for $n=2$, therefore we omit the details and list the resulting QLCs and their curvatures. The maps $\sigma$ although computed in the analysis are not listed as they are uniquely determined by $\nabla$ and the commutation relations. There is no case A as this did not have any quantum metrics. 

\medskip 
Case B (for metric $g_{B}$):\\
(B.1) 
$\nabla \extd e=0$,\quad $\nabla \extd x=0$,\quad
$\nabla \extd y= \mathrm{d} e\otimes \mathrm{d} x+ \mathrm{d}
x\otimes \mathrm{d} e+ \mathrm{d} e\otimes \mathrm{d} e+ \mathrm{d} x\otimes 
\mathrm{d} x$,\quad $R_\nabla=0$;\\
(B.2)
$\nabla \extd e=0$,\quad$\nabla \extd x= \mathrm{d} e\otimes \mathrm{d} y+ \mathrm{d}
y\otimes \mathrm{d} e+ \mathrm{d} e\otimes \mathrm{d} e+ \mathrm{d} y\otimes 
\mathrm{d} y$,\quad $\nabla \extd y=0$,\quad $R_{\nabla }=0$;\\
(B.3)
$\nabla \extd e=0$,\quad$\nabla \extd x= \mathrm{d} x\otimes \mathrm{d} y+ \mathrm{d}
y\otimes \mathrm{d} x+ \mathrm{d} x\otimes \mathrm{d} x+ \mathrm{d} y\otimes 
\mathrm{d} y=\nabla \extd y$,\quad $R_{\nabla }=0$.\\

Case C \newline
- for metric $g_{C.I}$:\newline
(C.I.1) $\nabla \mathrm{d} e =0, \quad \nabla \mathrm{d} x =0, \quad \nabla 
\mathrm{d} y = \mathrm{d} e\otimes \mathrm{d} e+\mathrm{d} e\otimes \mathrm{d
} x+\mathrm{d} x\otimes \mathrm{d} e+\mathrm{d} x\otimes \mathrm{d} x$;\quad $R_\nabla=0$;

(C.I.2) $\nabla \mathrm{d} e = \mathrm{d} e\otimes \mathrm{d} y+\mathrm{d}
y\otimes \mathrm{d} e+\mathrm{d} x\otimes \mathrm{d} y+\mathrm{d} y\otimes 
\mathrm{d} x, \quad \nabla \mathrm{d} x =0,$\newline
$\nabla \mathrm{d} y = \mathrm{d} e\otimes \mathrm{d} e+\mathrm{d} e\otimes 
\mathrm{d} x+\mathrm{d} x\otimes \mathrm{d} e+\mathrm{d} x\otimes \mathrm{d}
x+\mathrm{d} y\otimes \mathrm{d} y,\quad R_\nabla=0$;

%
%

(C.I.3) $\nabla \mathrm{d}e=\mathrm{d}e\otimes \mathrm{d}y+\mathrm{d}
y\otimes \mathrm{d}e,\quad \nabla \mathrm{d}x=\mathrm{d}y\otimes \mathrm{d}
y=\nabla \mathrm{d}y,\quad R_\nabla=0$; 
%
%
%
%

(C.I.4) $\nabla \mathrm{d} e = \mathrm{d} e\otimes \mathrm{d} y+\mathrm{d}
y\otimes \mathrm{d} e+\mathrm{d} x\otimes \mathrm{d} y+\mathrm{d} y\otimes 
\mathrm{d} x, \quad \nabla \mathrm{d} x =0, \newline
\nabla \mathrm{d} y = \mathrm{d} y\otimes \mathrm{d} y,\quad R_\nabla=0$; 
%
%
%
%

(C.I.5) $\nabla \mathrm{d} e = \mathrm{d} x\otimes \mathrm{d} y+\mathrm{d}
y\otimes \mathrm{d} x, \quad \nabla \mathrm{d} x = \mathrm{d} y\otimes 
\mathrm{d} y, \quad \nabla \mathrm{d} y =0,\quad R_\nabla=0$; 
%
%
%

(C.I.6) $\nabla \mathrm{d} e= \mathrm{d} y\otimes \mathrm{d} y, \quad \nabla 
\mathrm{d} x =0, \quad \nabla \mathrm{d} y =0,\quad R_\nabla=0$; 
%
%
%
%

(C.I.7) $\nabla \mathrm{d}e=\mathrm{d}y\otimes \mathrm{d}y+\mathrm{d}
e\otimes \mathrm{d}y+\mathrm{d}y\otimes \mathrm{d}e,\quad \nabla \mathrm{d}x=
\mathrm{d}y\otimes \mathrm{d}y=\nabla \mathrm{d}y,\quad R_\nabla=0$; 
%
%
%

(C.I.8) $\nabla \mathrm{d} e= \mathrm{d} y\otimes \mathrm{d} y+\mathrm{d}
e\otimes \mathrm{d} y+\mathrm{d} y\otimes \mathrm{d} e+\mathrm{d} x\otimes 
\mathrm{d} y+\mathrm{d} y\otimes \mathrm{d} x, \quad \nabla \mathrm{d} x=0,$ 
\newline
$\nabla \mathrm{d} y = \mathrm{d} y\otimes \mathrm{d} y,\quad R_\nabla=0$; 
%
%
%
%

(C.I.9) $\nabla \mathrm{d} e= \mathrm{d} y\otimes \mathrm{d} y+\mathrm{d}
x\otimes \mathrm{d} y+\mathrm{d} y\otimes \mathrm{d} x, \quad \nabla \mathrm{
d} x = \mathrm{d} y\otimes \mathrm{d} y, \quad \nabla \mathrm{d} y=0,\quad R_\nabla=0$; 
%
%
%

(C.I.10) $\nabla \mathrm{d}e=\mathrm{d}x\otimes \mathrm{d}x+\mathrm{d}
e\otimes \mathrm{d}y+\mathrm{d}y\otimes \mathrm{d}e,\quad \nabla \mathrm{d}x=
\mathrm{d}x\otimes \mathrm{d}y+\mathrm{d}y\otimes \mathrm{d}x$, \newline
$\nabla \mathrm{d}y=\mathrm{d}y\otimes \mathrm{d}y,\quad R_\nabla=0$; 
%
%
%

(C.I.11) $\nabla \mathrm{d}e=\mathrm{d}x\otimes \mathrm{d}x+\mathrm{d}
x\otimes \mathrm{d}y+\mathrm{d}y\otimes \mathrm{d}x,\quad \nabla \mathrm{d}x=
\mathrm{d}x\otimes \mathrm{d}y+\mathrm{d}y\otimes \mathrm{d}x,\quad \nabla \mathrm{d}y=0$,
\[ R_{\nabla }\mathrm{d}e=\mathrm{d}y\wedge \left( \mathrm{d}x\otimes \mathrm{d
}x\right) +\mathrm{d}x\wedge \left( \mathrm{d}y\otimes \mathrm{d}y\right)
,\quad R_{\nabla }\mathrm{d}x=\mathrm{d}x\wedge \left( \mathrm{d}y\otimes 
\mathrm{d}y\right) ,\quad R_{\nabla }\mathrm{d}y=0;\]

%
%
%
%

(C.I.12) $\nabla \mathrm{d}e=\mathrm{d}x\otimes \mathrm{d}x+\mathrm{d}
y\otimes \mathrm{d}y+\mathrm{d}e\otimes \mathrm{d}y+\mathrm{d}y\otimes 
\mathrm{d}e,\quad \nabla \mathrm{d}x=\mathrm{d}x\otimes \mathrm{d}y+\mathrm{d
}y\otimes \mathrm{d}x$,  
$\nabla \mathrm{d}y=\mathrm{d}y\otimes \mathrm{d}y,\quad R_\nabla=0$; 
%
%
%

(C.I.13) $\nabla \mathrm{d}e=\mathrm{d}x\otimes \mathrm{d}x+\mathrm{d}
y\otimes \mathrm{d}y+\mathrm{d}x\otimes \mathrm{d}y+\mathrm{d}y\otimes 
\mathrm{d}x,\quad \nabla \mathrm{d}x=\mathrm{d}x\otimes \mathrm{d}y+\mathrm{d
}y\otimes \mathrm{d}x,\quad \nabla \mathrm{d}y=0$,
\[ R_{\nabla }\mathrm{d}e=\mathrm{d}x\wedge \left( \mathrm{d}y\otimes \mathrm{d
}y\right) +\mathrm{d}y\otimes (\mathrm{d}x\otimes \mathrm{d}x),\quad R_{\nabla }\mathrm{d}x=\mathrm{
d}x\wedge \left( \mathrm{d}y\otimes \mathrm{d}y\right) ,\quad R_{\nabla }
\mathrm{d}y=0.\] 
%
%
%
%
%

- for metric $g_{C.II}$:\newline
(C.II.1) $\nabla \mathrm{d} e=0$,\quad $\nabla \mathrm{d} x= \mathrm{d}
e\otimes \mathrm{d} x+ \mathrm{d} x\otimes \mathrm{d} e+ \mathrm{d} e\otimes 
\mathrm{d} y+ \mathrm{d} y\otimes \mathrm{d} e$,\newline
$\nabla \mathrm{d} y= \mathrm{d} e\otimes \mathrm{d} x+ \mathrm{d} x\otimes 
\mathrm{d} e+ \mathrm{d} e\otimes \mathrm{d} y+ \mathrm{d} y\otimes \mathrm{d
} e+ \mathrm{d} e\otimes \mathrm{d} e,\quad R_\nabla=0$;

(C.II.2) $\nabla \mathrm{d}e=\mathrm{d}e\otimes \mathrm{d}y+\mathrm{d}
y\otimes \mathrm{d}e+\mathrm{d}x\otimes \mathrm{d}y+\mathrm{d}y\otimes 
\mathrm{d}x$,\quad $\nabla \mathrm{d}x=\mathrm{d}e\otimes \mathrm{d}x+
\mathrm{d}x\otimes \mathrm{d}e+\mathrm{d}e\otimes \mathrm{d}y+\mathrm{d}
y\otimes \mathrm{d}e$,\quad $\nabla \mathrm{d}y=\mathrm{d}e\otimes \mathrm{d}
x+\mathrm{d}x\otimes \mathrm{d}e+\mathrm{d}e\otimes \mathrm{d}y+\mathrm{d}
y\otimes \mathrm{d}e+\mathrm{d}e\otimes \mathrm{d}e+\mathrm{d}y\otimes 
\mathrm{d}y$,
\[ R_{\nabla }\mathrm{d}e=\mathrm{d}y\wedge \left( \mathrm{d}e\otimes \mathrm{d
}e\right) +\mathrm{d}e\wedge \left( \mathrm{d}x\otimes \mathrm{d}x+\mathrm{d}
y\otimes \mathrm{d}y+\mathrm{d}x\otimes \mathrm{d}y+\mathrm{d}y\otimes 
\mathrm{d}x\right),\]
\[R_{\nabla }\mathrm{\mathrm{d}}x=\mathrm{d}y\wedge \left( \mathrm{d}x\otimes 
\mathrm{d}x\right) +\mathrm{d}x\wedge \left( \mathrm{d}e\otimes \mathrm{d}y+
\mathrm{d}y\otimes \mathrm{d}e+\mathrm{d}y\otimes \mathrm{d}y\right) ,\]
\[ R_{\nabla }\mathrm{\mathrm{d}}y=\mathrm{d}e\wedge \left( \mathrm{d}y\otimes 
\mathrm{d}y\right) +\mathrm{d}x\wedge \left( \mathrm{d}y\otimes \mathrm{d}
y\right) +\mathrm{d}y\wedge \left( \mathrm{d}x\otimes \mathrm{d}x\right).\] 

The remaining solutions (C.II.3) - (C.II.13) are the same as for the
metric $g_{C.I}$, i.e. are equal to cases (C.I.3) - (C.I.13) respectively.
\newline

Case D\newline
- for metric $g_{D.I}$:\newline
(D.I.1) $\nabla \mathrm{d} e=0$,\quad $\nabla \mathrm{d} x= \mathrm{d}
e\otimes \mathrm{d} e $,\quad $\nabla \mathrm{d} y= \mathrm{d} e\otimes 
\mathrm{d} x+ \mathrm{d} x\otimes \mathrm{d} e,\quad R_\nabla=0$;

(D.I.2) $\nabla \mathrm{d} e= \mathrm{d} x\otimes \mathrm{d} y+ \mathrm{d}
y\otimes \mathrm{d} x+ \mathrm{d} y\otimes \mathrm{d} y $,\quad$\nabla 
\mathrm{d} x= \mathrm{d} x\otimes \mathrm{d} y+ \mathrm{d} y\otimes \mathrm{d
} x $,\quad$\nabla \mathrm{d} y= \mathrm{d} y\otimes \mathrm{d} y,\quad R_\nabla=0$;

(D.I.3) $\nabla \mathrm{d}e=\mathrm{d}e\otimes \mathrm{d}x+\mathrm{d}
x\otimes \mathrm{d}e+\mathrm{d}e\otimes \mathrm{d}y+\mathrm{d}y\otimes 
\mathrm{d}e+\mathrm{d}e\otimes \mathrm{d}e+\mathrm{d}y\otimes \mathrm{d}y$,
\newline
$\nabla \mathrm{d}x=\mathrm{d}e\otimes \mathrm{d}x+\mathrm{d}x\otimes 
\mathrm{d}e+\mathrm{d}e\otimes \mathrm{d}e+\mathrm{d}y\otimes \mathrm{d}y$
,\\ $\nabla \mathrm{d}y=\mathrm{d}e\otimes \mathrm{d}e+\mathrm{d}x\otimes 
\mathrm{d}y+\mathrm{d}y\otimes \mathrm{d}x+\mathrm{d}y\otimes \mathrm{d}y$,
\[ R_{\nabla }\mathrm{d}e=\mathrm{d}x\wedge \left( \mathrm{d}e\otimes \mathrm{d
}e\right) +\mathrm{d}e\wedge \left( \mathrm{d}x\otimes \mathrm{d}x\right) +
\mathrm{d}y\wedge \mathrm{d}e\otimes \mathrm{d}e,\]
\[R_{\nabla }\mathrm{d}x=\mathrm{d}e\wedge \left( \mathrm{d}x\otimes \mathrm{d
}x\right) +\mathrm{d}x\wedge \left( \mathrm{d}e\otimes \mathrm{d}e+\mathrm{d}
y\otimes \mathrm{d}e+\mathrm{d}e\otimes \mathrm{d}y\right),\]
\[ R_{\nabla }\mathrm{d}y=\mathrm{d}e\wedge (\mathrm{d}y\otimes \mathrm{d}y)+
\mathrm{d}y\wedge \left( \mathrm{d}e\otimes \mathrm{d}e+\mathrm{d}e\otimes 
\mathrm{d}x+\mathrm{d}x\otimes \mathrm{d}x+\mathrm{d}x\otimes \mathrm{d}
e\right).\]

- metric $g_{D.II}$:\newline
(D.II.1) $\nabla \mathrm{d} e=0 $,\quad$\nabla \mathrm{d} x= \mathrm{d}
e\otimes \mathrm{d} x+ \mathrm{d} x\otimes \mathrm{d} e+ \mathrm{d} e\otimes 
\mathrm{d} y+ \mathrm{d} y\otimes \mathrm{d} e+ \mathrm{d} x\otimes \mathrm{d
} y+ \mathrm{d} y\otimes \mathrm{d} x $,\newline
$\nabla \mathrm{d} y= \mathrm{d} e\otimes \mathrm{d} y+ \mathrm{d} y\otimes 
\mathrm{d} e+ \mathrm{d} e\otimes \mathrm{d} e+ \mathrm{d} y\otimes \mathrm{d
} y,\quad R_\nabla=0$;

(D.II.2) $\nabla \mathrm{d} e=0 $,\quad$\nabla \mathrm{d} x= \mathrm{d}
e\otimes \mathrm{d} x+ \mathrm{d} x\otimes \mathrm{d} e+ \mathrm{d} e\otimes 
\mathrm{d} e+ \mathrm{d} x\otimes \mathrm{d} x $,\newline
$\nabla \mathrm{d} y= \mathrm{d} e\otimes \mathrm{d} x+ \mathrm{d} x\otimes 
\mathrm{d} e+ \mathrm{d} e\otimes \mathrm{d} y+ \mathrm{d} y\otimes \mathrm{d
} e+ \mathrm{d} x\otimes \mathrm{d} y+ \mathrm{d} y\otimes \mathrm{d} x,\quad R_\nabla=0$;

(D.II.3) $\nabla \mathrm{d} e=0 $,\quad$\nabla \mathrm{d} x= \mathrm{d}
x\otimes \mathrm{d} y+ \mathrm{d} y\otimes \mathrm{d} x+ \mathrm{d} x\otimes 
\mathrm{d} x+ \mathrm{d} y\otimes \mathrm{d} y $,\newline
$\nabla \mathrm{d} y= \mathrm{d} x\otimes \mathrm{d} y+ \mathrm{d} y\otimes 
\mathrm{d} x+ \mathrm{d} x\otimes \mathrm{d} x+ \mathrm{d} y\otimes \mathrm{d
} y,\quad R_\nabla=0$.

- for metric $g_{D.III}$\newline
(D.III.1) $\nabla \mathrm{d} e=0$,\quad $\nabla \mathrm{d} x= \mathrm{d}
e\otimes \mathrm{d} y + \mathrm{d} y\otimes \mathrm{d} e $,\quad $\nabla 
\mathrm{d} y= \mathrm{d} e\otimes \mathrm{d} e,\quad R_\nabla=0$;

(D.III.2) $\nabla \mathrm{d} e= \mathrm{d} x\otimes \mathrm{d} x + \mathrm{d}
x\otimes \mathrm{d} y + \mathrm{d} y\otimes \mathrm{d} x$,\quad $\nabla 
\mathrm{d} x= \mathrm{d} x\otimes \mathrm{d} x $,\newline
$\nabla \mathrm{d} y= \mathrm{d} x\otimes \mathrm{d} y + \mathrm{d} y\otimes 
\mathrm{d} x,\quad R_\nabla=0$;

(D.III.3) $\nabla \mathrm{d}e=\mathrm{d}e\otimes \mathrm{d}e+\mathrm{d}
e\otimes \mathrm{d}x+\mathrm{d}e\otimes \mathrm{d}y+\mathrm{d}x\otimes 
\mathrm{d}e+\mathrm{d}x\otimes \mathrm{d}x+\mathrm{d}y\otimes \mathrm{d}e$
,\\ $\nabla \mathrm{d}x=\mathrm{d}e\otimes \mathrm{d}e+\mathrm{d}x\otimes 
\mathrm{d}x+\mathrm{d}x\otimes \mathrm{d}y+\mathrm{d}y\otimes \mathrm{d}x$,
\quad $\nabla \mathrm{d}y=\mathrm{d}e\otimes \mathrm{d}e+\mathrm{d}e\otimes 
\mathrm{d}y+\mathrm{d}x\otimes \mathrm{d}x+\mathrm{d}y\otimes \mathrm{d}e$,
\[ R_{\nabla }\mathrm{d}e=\left( \mathrm{d}x+\mathrm{d}y\right) \wedge \left( 
\mathrm{d}e\otimes \mathrm{d}e\right) +\mathrm{d}e\wedge \left( \mathrm{d}
y\otimes \mathrm{d}y\right),\]
\[ R_{\nabla }\mathrm{d}x=\mathrm{d}e\wedge \left( \mathrm{d}x\otimes \mathrm{d
}x\right) +\mathrm{d}x\wedge \left( \mathrm{d}e\otimes \mathrm{d}e+\mathrm{d}
y\otimes \mathrm{d}y+\mathrm{d}e\otimes \mathrm{d}y+\mathrm{d}y\otimes 
\mathrm{d}e\right),\]
\[ R_{\nabla }\mathrm{d}y=\mathrm{d}e\wedge \left( \mathrm{d}y\otimes 
\mathrm{d}y\right) +\mathrm{d}y\wedge \left( \mathrm{d}e\otimes \mathrm{d}e+
\mathrm{d}x\otimes \mathrm{d}e+\mathrm{d}e\otimes \mathrm{d}x\right).\]

Case E\newline
- for metric $g_{E.I}$:\newline
(E.I.1) $\nabla \mathrm{d}e=0$,\quad $\nabla \mathrm{d}x=\mathrm{d}e\otimes 
\mathrm{d}x+\mathrm{d}x\otimes \mathrm{d}e+\mathrm{d}e\otimes \mathrm{d}y+
\mathrm{d}y\otimes \mathrm{d}e$,\newline
$\nabla \mathrm{d}y=\mathrm{d}e\otimes \mathrm{d}x+\mathrm{d}x\otimes 
\mathrm{d}e+\mathrm{d}e\otimes \mathrm{d}y+\mathrm{d}y\otimes \mathrm{d}e+
\mathrm{d}e\otimes \mathrm{d}e,\quad R_\nabla=0;$

(E.I.2) $\nabla \mathrm{d}e=\mathrm{d}e\otimes \mathrm{d}y+\mathrm{d}
y\otimes \mathrm{d}e$,\quad $\nabla \mathrm{d}x=\mathrm{d}y\otimes \mathrm{d}
y=\nabla \mathrm{d}y,\quad R_\nabla=0$;

(E.I.3) $\nabla \mathrm{d} e= \mathrm{d} e\otimes \mathrm{d} y+ \mathrm{d}
y\otimes \mathrm{d} e+ \mathrm{d} x\otimes \mathrm{d} y+ \mathrm{d} y\otimes 
\mathrm{d} x $,\quad $\nabla \mathrm{d} x=0 $,\quad $\nabla \mathrm{d} y= 
\mathrm{d} y\otimes \mathrm{d} y,\quad R_\nabla=0$;

(E.I.4) $\nabla \mathrm{d}e=\mathrm{d}e\otimes \mathrm{d}y+\mathrm{d}
y\otimes \mathrm{d}e+\mathrm{d}x\otimes \mathrm{d}y+\mathrm{d}y\otimes 
\mathrm{d}x$,\quad $\nabla \mathrm{d}x=\mathrm{d}e\otimes \mathrm{d}x+
\mathrm{d}x\otimes \mathrm{d}e+\mathrm{d}e\otimes \mathrm{d}y+\mathrm{d}
y\otimes \mathrm{d}e$,\quad $\nabla \mathrm{d}y=\mathrm{d}e\otimes \mathrm{d}
x+\mathrm{d}x\otimes \mathrm{d}e+\mathrm{d}e\otimes \mathrm{d}y+\mathrm{d}
y\otimes \mathrm{d}e+\mathrm{d}e\otimes \mathrm{d}e+\mathrm{d}y\otimes 
\mathrm{d}y$,
\[ R_{\nabla }\mathrm{d}e=\mathrm{d}y\wedge \left( \mathrm{d}e\otimes \mathrm{d
}e\right) +\mathrm{d}e\wedge \left( \mathrm{d}x\otimes \mathrm{d}x+\mathrm{d}
x\otimes \mathrm{d}y+\mathrm{d}y\otimes \mathrm{d}x+\mathrm{d}y\otimes 
\mathrm{d}y\right),\]
\[R_{\nabla }\mathrm{d}x=\mathrm{d}x\wedge (\mathrm{d}e\otimes \mathrm{d}y+
\mathrm{d}y\otimes \mathrm{d}e)+\mathrm{d}x\wedge \mathrm{d}y\otimes \mathrm{
d}y+\mathrm{d}y\wedge \mathrm{d}x\otimes \mathrm{d}x,\]
\[  R_{\nabla }\mathrm{d}y=\mathrm{d}e\wedge \left( \mathrm{d}y\otimes \mathrm{d
}y\right) +\mathrm{d}y\wedge \left( \mathrm{d}x\otimes \mathrm{d}x\right) +
\mathrm{d}x\wedge \left( \mathrm{d}y\otimes \mathrm{d}y\right);\]

(E.I.5) $\nabla \mathrm{d} e= \mathrm{d} x\otimes \mathrm{d} y+ \mathrm{d}
y\otimes \mathrm{d} x $ ,\quad $\nabla \mathrm{d} x= \mathrm{d} y\otimes 
\mathrm{d} y $,\quad $\nabla \mathrm{d} y=0,\quad R_\nabla=0;$

(E.I.6) $\nabla \mathrm{d} e= \mathrm{d} y\otimes \mathrm{d} y $,\quad $
\nabla \mathrm{d} x=0 $ ,\quad $\nabla \mathrm{d} y=0,\quad R_\nabla=0$;

(E.I.7) $\nabla \mathrm{d} e= \mathrm{d} e\otimes \mathrm{d} y + \mathrm{d}
y\otimes \mathrm{d} e + \mathrm{d} y\otimes \mathrm{d} y $,\quad $\nabla 
\mathrm{d} x= \mathrm{d} y\otimes \mathrm{d} y $ ,\quad $\nabla \mathrm{d}
y= \mathrm{d} y\otimes \mathrm{d} y,\quad R_\nabla=0$;

(E.I.8) $\nabla \mathrm{d} e= \mathrm{d} e\otimes \mathrm{d} y + \mathrm{d}
x\otimes \mathrm{d} y + \mathrm{d} y\otimes \mathrm{d} e + \mathrm{d}
y\otimes \mathrm{d} x + \mathrm{d} y\otimes \mathrm{d} y $ ,\quad $\nabla 
\mathrm{d} x=0 $ ,\quad $\nabla \mathrm{d} y= \mathrm{d} y\otimes \mathrm{d}
y,\quad R_\nabla=0$;

(E.I.9) $\nabla \mathrm{d} e= \mathrm{d} x\otimes \mathrm{d} y + \mathrm{d}
y\otimes \mathrm{d} x + \mathrm{d} y\otimes \mathrm{d} y $,\quad $\nabla 
\mathrm{d} x= \mathrm{d} y\otimes \mathrm{d} y $,\quad $\nabla \mathrm{d}
y=0,\quad R_\nabla=0$;

(E.I.10) $\nabla \mathrm{d} e= \mathrm{d} e\otimes \mathrm{d} y + \mathrm{d}
x\otimes \mathrm{d} x + \mathrm{d} y\otimes \mathrm{d} e $ ,\quad $\nabla 
\mathrm{d} x= \mathrm{d} x\otimes \mathrm{d} y + \mathrm{d} y\otimes \mathrm{%
d} x $ ,\newline
$\nabla \mathrm{d} y= \mathrm{d} y\otimes \mathrm{d} y,\quad R_\nabla=0$;

(E.I.11) $\nabla \mathrm{d} e= \mathrm{d} x\otimes \mathrm{d} x + \mathrm{d}
x\otimes \mathrm{d} y + \mathrm{d} y\otimes \mathrm{d} x $ ,\quad $\nabla 
\mathrm{d} x= \mathrm{d} x\otimes \mathrm{d} y + \mathrm{d} y\otimes \mathrm{%
d} x $ ,\quad $\nabla \mathrm{d} y=0,\quad R_\nabla=0$;

(E.I.12) $\nabla \mathrm{d} e= \mathrm{d} e\otimes \mathrm{d} y + \mathrm{d}
x\otimes \mathrm{d} x + \mathrm{d} y\otimes \mathrm{d} e + \mathrm{d}
y\otimes \mathrm{d} y $ ,\quad $\nabla \mathrm{d} x= \mathrm{d} x\otimes 
\mathrm{d} y + \mathrm{d} y\otimes \mathrm{d} x $,\newline
$\nabla \mathrm{d} y= \mathrm{d} y\otimes \mathrm{d} y,\quad R_\nabla=0$;

(E.I.13) $\nabla \mathrm{d}e=\mathrm{d}x\otimes \mathrm{d}x+\mathrm{d}
x\otimes \mathrm{d}y+\mathrm{d}y\otimes \mathrm{d}x+\mathrm{d}y\otimes 
\mathrm{d}y$ ,\quad $\nabla \mathrm{d}x=\mathrm{d}x\otimes \mathrm{d}y+
\mathrm{d}y\otimes \mathrm{d}x$,
$\nabla \mathrm{d}y=0$,
\[ R_{\nabla }\mathrm{d}e=\mathrm{d}x\wedge \left( \mathrm{d}y\otimes \mathrm{d
}y\right) +\mathrm{d}y\wedge \left( \mathrm{d}x\otimes \mathrm{d}x\right)
,\quad R_{\nabla }\mathrm{d}x=\mathrm{d}x\wedge \left( \mathrm{d}y\otimes 
\mathrm{d}y\right) ,\quad R_{\nabla }\mathrm{d}y=0.\]

- for metric $g_{E.II}$:\newline
(E.II.1) $\nabla \mathrm{d} e=0 $ ,\quad $\nabla \mathrm{d} x= \mathrm{d} e
\otimes \mathrm{d} e + \mathrm{d} e \otimes \mathrm{d} x + \mathrm{d} x
\otimes \mathrm{d} e + \mathrm{d} x \otimes \mathrm{d} x $,\newline
$\nabla \mathrm{d} y= \mathrm{d} e \otimes \mathrm{d} x + \mathrm{d} e
\otimes \mathrm{d} y + \mathrm{d} x \otimes \mathrm{d} e + \mathrm{d} x
\otimes \mathrm{d} y + \mathrm{d} y \otimes \mathrm{d} e + \mathrm{d} y
\otimes \mathrm{d} x,\quad R_\nabla=0$;

(E.II.2) $\nabla \mathrm{d}e=\mathrm{d}x\otimes \mathrm{d}x$,\quad $\nabla 
\mathrm{d}x=\mathrm{d}x\otimes \mathrm{d}y+\mathrm{d}y\otimes \mathrm{d}x+
\mathrm{d}y\otimes \mathrm{d}y$,\quad $\nabla \mathrm{d}y=0$,
\[ R_{\nabla }\mathrm{d}e=\mathrm{d}y\wedge \left( \mathrm{d}x\otimes \mathrm{d
}x\right) +\mathrm{d}x\wedge \left( \mathrm{d}y\otimes \mathrm{d}y\right)
,\quad R_{\nabla }\mathrm{d}x=\mathrm{d}x\wedge \left( \mathrm{d}y\otimes 
\mathrm{d}y\right),\quad R_{\nabla }\mathrm{d}y=0;\]

(E.II.3) $\nabla \mathrm{d}e=\mathrm{d}x\otimes \mathrm{d}y+\mathrm{d}
y\otimes \mathrm{d}x$,\quad $\nabla \mathrm{d}x=\mathrm{d}e\otimes \mathrm{d}
e+\mathrm{d}e\otimes \mathrm{d}x+\mathrm{d}x\otimes \mathrm{d}e+\mathrm{d}
x\otimes \mathrm{d}x+\mathrm{d}y\otimes \mathrm{d}y$,\newline
$\nabla \mathrm{d}y=\mathrm{d}e\otimes \mathrm{d}x+\mathrm{d}e\otimes 
\mathrm{d}y+\mathrm{d}x\otimes \mathrm{d}e+\mathrm{d}x\otimes \mathrm{d}y+
\mathrm{d}y\otimes \mathrm{d}e+\mathrm{d}y\otimes \mathrm{d}x$,
\[ R_{\nabla }\mathrm{d}e=\mathrm{d}y\wedge \left( \mathrm{d}e\otimes \mathrm{d
}e\right) +\mathrm{d}e\wedge \left( \mathrm{d}x\otimes \mathrm{d}x+\mathrm{d}
y\otimes \mathrm{d}x+\mathrm{d}x\otimes \mathrm{d}y\right),\]
\[ R_{\nabla }
\mathrm{d}x=\mathrm{d}x\wedge \left( \mathrm{d}e\otimes \mathrm{d}y\right) +
\mathrm{d}y\wedge \left( \mathrm{d}x\otimes \mathrm{d}x\right) +\mathrm{d}
x\wedge \left( \mathrm{d}y\otimes \mathrm{d}e\right),\]
\[ R_{\nabla }\mathrm{d}y=\mathrm{d}y\wedge \left( \mathrm{d}x\otimes \mathrm{d%
}x\right) +\mathrm{d}x\wedge \left( \mathrm{d}y\otimes \mathrm{d}y\right) +
\mathrm{d}e\wedge \left( \mathrm{d}y\otimes \mathrm{d}y\right);\]

(E.II.4) $\nabla \mathrm{d}e=\mathrm{d}x\otimes \mathrm{d}x+\mathrm{d}%
y\otimes \mathrm{d}y$,\quad $\nabla \mathrm{d}x=\mathrm{d}x\otimes \mathrm{d}%
y+\mathrm{d}y\otimes \mathrm{d}x+\mathrm{d}y\otimes \mathrm{d}y$,\quad $%
\nabla \mathrm{d}y=0,$
\[ R_{\nabla }\mathrm{d}e=\mathrm{d}y\wedge \left( \mathrm{d}x\otimes \mathrm{d%
}x\right) +\mathrm{d}x\wedge \left( \mathrm{d}y\otimes \mathrm{d}y\right)
,\quad R_{\nabla }\mathrm{d}x=\mathrm{d}x\wedge \left( \mathrm{d}y\otimes 
\mathrm{d}y\right),\quad R_{\nabla }\mathrm{d}y=0.\]

And the remaining ones are: (E.II.5.)=(E.I.11.), (E.II.6.)=(E.I.5.),
(E.II.7.)=(E.I.2.), (E.II.8.)=(E.I.3.), (E.II.9.)=(E.I.6.),
(E.II.10.)=(E.I.13.), (E.II.11.)=(E.I.9.), (E.II.12.)=(E.I.7.),
(E.II.13.)=(E.I.8.).

- for metric $g_{E.III}$

(E.III.1) $\nabla \mathrm{d}e=\mathrm{d}x\otimes \mathrm{d}x$,\quad $\nabla 
\mathrm{d}x=\mathrm{d}x\otimes \mathrm{d}y+\mathrm{d}y\otimes \mathrm{d}x$
,\quad $\nabla \mathrm{d}y=0,$
\[ R_{\nabla }\mathrm{d}e=\mathrm{d}y\wedge \left( \mathrm{d}x\otimes \mathrm{d
}x\right) ,\quad R_{\nabla }\mathrm{d}x=\mathrm{d}x\wedge \left( \mathrm{d}
y\otimes \mathrm{d}y\right) ,\quad R_{\nabla }\mathrm{d}y=0;\]

(E.III.2) $\nabla \mathrm{d}e=0$,\quad $\nabla \mathrm{d}x=\mathrm{d}
e\otimes \mathrm{d}e+\mathrm{d}e\otimes \mathrm{d}y+\mathrm{d}x\otimes 
\mathrm{d}x+\mathrm{d}y\otimes \mathrm{d}e+\mathrm{d}y\otimes \mathrm{d}y$,
\newline
$\nabla \mathrm{d}y=\mathrm{d}e\otimes \mathrm{d}x+\mathrm{d}x\otimes 
\mathrm{d}e+\mathrm{d}x\otimes \mathrm{d}y+\mathrm{d}y\otimes \mathrm{d}x,\quad R_\nabla=0$;

(E.III.3) $\nabla \mathrm{d}e=\mathrm{d}x\otimes \mathrm{d}y+\mathrm{d}
y\otimes \mathrm{d}x$,\quad $\nabla \mathrm{d}x=\mathrm{d}e\otimes \mathrm{d}
e+\mathrm{d}e\otimes \mathrm{d}y+\mathrm{d}x\otimes \mathrm{d}x+\mathrm{d}
y\otimes \mathrm{d}e$,\quad $\nabla \mathrm{d}y=\mathrm{d}e\otimes \mathrm{d}
x+\mathrm{d}x\otimes \mathrm{d}e+\mathrm{d}x\otimes \mathrm{d}y+\mathrm{d}
y\otimes \mathrm{d}x,$
\[ R_{\nabla }\mathrm{d}e=\mathrm{d}y\wedge \left( \mathrm{d}e\otimes \mathrm{d
}e\right) +\mathrm{d}e\wedge \left( \mathrm{d}x\otimes \mathrm{d}x+\mathrm{d}
y\otimes \mathrm{d}y\right),\quad R_{\nabla }\mathrm{d}x=\mathrm{d}x\wedge \left( \mathrm{d}e\otimes \mathrm{d
}y+\mathrm{d}y\otimes \mathrm{d}e+\mathrm{d}y\otimes \mathrm{d}y\right),\]
\[ R_{\nabla }\mathrm{d}y=\mathrm{d}e\wedge \left( \mathrm{d}y\otimes \mathrm{d
}y\right) +\mathrm{d}y\wedge \left( \mathrm{d}y\otimes \mathrm{d}y\right);\]

(E.III.4) $\nabla \mathrm{d}e=\mathrm{d}x\otimes \mathrm{d}x+\mathrm{d}
x\otimes \mathrm{d}y+\mathrm{d}y\otimes \mathrm{d}x$,\quad $\nabla \mathrm{d}
x=\mathrm{d}x\otimes \mathrm{d}y+\mathrm{d}y\otimes \mathrm{d}x+\mathrm{d}
y\otimes \mathrm{d}y$,\quad $\nabla \mathrm{d}y=0$,
\[ R_{\nabla }\mathrm{d}e=\mathrm{d}y\wedge \left( \mathrm{d}x\otimes \mathrm{d
}x\right) ,\quad R_{\nabla }\mathrm{d}x=\mathrm{d}x\wedge \left( \mathrm{d}
y\otimes \mathrm{d}y\right) ,\quad R_{\nabla }\mathrm{d}y=0;\]

(E.III.5) $\nabla \mathrm{d}e=\mathrm{d}e\otimes \mathrm{d}y+\mathrm{d}
y\otimes \mathrm{d}e$,\quad $\nabla \mathrm{d}x=0$,\quad $\nabla \mathrm{d}y=
\mathrm{d}y\otimes \mathrm{d}y,\quad R_\nabla=0$;

(E.III.6) $\nabla \mathrm{d}e=\mathrm{d}e\otimes \mathrm{d}y+\mathrm{d}%
x\otimes \mathrm{d}y+\mathrm{d}y\otimes \mathrm{d}e+\mathrm{d}y\otimes 
\mathrm{d}x$,\quad $\nabla \mathrm{d}x=\mathrm{d}y\otimes \mathrm{d}y=\nabla 
\mathrm{d}y,\quad R_\nabla=0$;

(E.III.7) $\nabla \mathrm{d}e=\mathrm{d}x\otimes \mathrm{d}x+\mathrm{d}
y\otimes \mathrm{d}y$,\quad $\nabla \mathrm{d}x=\mathrm{d}x\otimes \mathrm{d}
y+\mathrm{d}y\otimes \mathrm{d}x$,\quad $\nabla \mathrm{d}y=0$,
\[ R_{\nabla }\mathrm{d}e=\mathrm{d}y\wedge \left( \mathrm{d}x\otimes \mathrm{d
}x\right) ,\quad R_{\nabla }\mathrm{d}x=\mathrm{d}x\wedge \left( \mathrm{d}
y\otimes \mathrm{d}y\right) ,\quad R_{\nabla }\mathrm{d}y=0;\]

(E.III.8) $\nabla \mathrm{d}e=\mathrm{d}x\otimes \mathrm{d}x+\mathrm{d}
x\otimes \mathrm{d}y+\mathrm{d}y\otimes \mathrm{d}x+\mathrm{d}y\otimes 
\mathrm{d}y$,\newline
$\nabla \mathrm{d}x=\mathrm{d}x\otimes \mathrm{d}y+\mathrm{d}y\otimes 
\mathrm{d}x+\mathrm{d}y\otimes \mathrm{d}y$,\quad $\nabla \mathrm{d}y=0,$
\[ R_{\nabla }\mathrm{d}e=\mathrm{d}y\wedge \left( \mathrm{d}x\otimes \mathrm{d
}x\right) ,\quad R_{\nabla }\mathrm{d}x=\mathrm{d}x\wedge \left( \mathrm{d}
y\otimes \mathrm{d}y\right) ,\quad R_{\nabla }\mathrm{d}y=0;\]

(E.III.9) $\nabla \mathrm{d}e=\mathrm{d}e\otimes \mathrm{d}y+\mathrm{d}
y\otimes \mathrm{d}e+\mathrm{d}y\otimes \mathrm{d}y$,\quad $\nabla \mathrm{d}
x=0$ ,\quad $\nabla \mathrm{d}y=\mathrm{d}y\otimes \mathrm{d}y,\quad R_\nabla=0$;

(E.III.10) $\nabla \mathrm{d}e=\mathrm{d}e\otimes \mathrm{d}y+\mathrm{d}%
x\otimes \mathrm{d}y+\mathrm{d}y\otimes \mathrm{d}e+\mathrm{d}y\otimes 
\mathrm{d}x+\mathrm{d}y\otimes \mathrm{d}y$,\quad $\nabla \mathrm{d}x=%
\mathrm{d}y\otimes \mathrm{d}y$,\quad $\nabla \mathrm{d}y=\mathrm{d}y\otimes 
\mathrm{d}y,\quad R_\nabla=0$.

The remaining ones are: (E.III.11) = (E.I.5), (E.III.12)=(E.I.6),
(E.III.13)=(E.I.9).

- for metric $g_{E.IV}$:

(E.IV.1) $\nabla \mathrm{d} e=0 $,\quad $\nabla \mathrm{d} x= \mathrm{d}
e\otimes \mathrm{d} e + \mathrm{d} e\otimes \mathrm{d} x + \mathrm{d}
e\otimes \mathrm{d} y + \mathrm{d} x\otimes \mathrm{d} e + \mathrm{d}
y\otimes \mathrm{d} e $,\newline
$\nabla \mathrm{d} y= \mathrm{d} e\otimes \mathrm{d} x + \mathrm{d} e\otimes 
\mathrm{d} y + \mathrm{d} x\otimes \mathrm{d} e + \mathrm{d} y\otimes 
\mathrm{d} e,\quad R_\nabla=0$;

%
%
%
%
%

(E.IV.2) $\nabla \mathrm{d}e=\mathrm{d}x\otimes \mathrm{d}y+\mathrm{d}%
y\otimes \mathrm{d}x$,\quad $\nabla \mathrm{d}x=\mathrm{d}e\otimes \mathrm{d}%
e+\mathrm{d}e\otimes \mathrm{d}x+\mathrm{d}e\otimes \mathrm{d}y+\mathrm{d}%
x\otimes \mathrm{d}e+\mathrm{d}y\otimes \mathrm{d}e+\mathrm{d}y\otimes 
\mathrm{d}y$,\quad $\nabla \mathrm{d}y=\mathrm{d}e\otimes \mathrm{d}x+%
\mathrm{d}e\otimes \mathrm{d}y+\mathrm{d}x\otimes \mathrm{d}e+\mathrm{d}%
y\otimes \mathrm{d}e,$
\[R_{\nabla }\mathrm{d}e=\mathrm{d}y\wedge \left( \mathrm{d}e\otimes \mathrm{d
}e\right) +\mathrm{d}e\wedge \left( \mathrm{d}x\otimes \mathrm{d}x+\mathrm{d}
x\otimes \mathrm{d}y+\mathrm{d}y\otimes \mathrm{d}x+\mathrm{d}y\otimes 
\mathrm{d}y\right),\]
\[ R_{\nabla }\mathrm{d}x=\mathrm{d}x\wedge \left( 
\mathrm{d}e\otimes \mathrm{d}y+\mathrm{d}y\otimes \mathrm{d}y\right) +
\mathrm{d}y\wedge \left( \mathrm{d}x\otimes \mathrm{d}x\right) +\mathrm{d}
x\wedge \left( \mathrm{d}y\otimes \mathrm{d}e\right),\]
\[  R_{\nabla }
\mathrm{d}y=\mathrm{d}e\wedge \left( \mathrm{d}y\otimes \mathrm{d}y\right) +
\mathrm{d}y\wedge \left( \mathrm{d}x\otimes \mathrm{d}x\right) +\mathrm{d}
x\wedge \left( \mathrm{d}y\otimes \mathrm{d}y\right);\]

And (E.IV.3)=(E.I.11), (E.IV.4)=(E.I.5),(E.IV.5)=(E.II.2), (E.IV.6)=(E.I.2),
(E.IV.7)=(E.I.3), (E.IV.8)=(E.I.6), (E.IV.9)=(E.II.4), (E.IV.10)=(E.I.13),
(E.IV.11)=(E.I.9), (E.IV.12)=(E.I.7), (E.IV.13)=(E.I.8).

Case F\newline
- for metric $g_{F.I}$

(F.I.1) $\nabla \mathrm{d} e=0 $,\quad $\nabla \mathrm{d} x= \mathrm{d}
e\otimes \mathrm{d} e $,\quad $\nabla \mathrm{d} y= \mathrm{d} e\otimes 
\mathrm{d} e + \mathrm{d} e\otimes \mathrm{d} x + \mathrm{d} x\otimes 
\mathrm{d} e,\quad R_\nabla=0$;

(F.I.2) $\nabla \mathrm{d}e=\mathrm{d}e\otimes \mathrm{d}x+\mathrm{d}
e\otimes \mathrm{d}y+\mathrm{d}x\otimes \mathrm{d}e+\mathrm{d}y\otimes 
\mathrm{d}e+\mathrm{d}y\otimes \mathrm{d}y$,\quad $\nabla \mathrm{d}x=
\mathrm{d}e\otimes \mathrm{d}y+\mathrm{d}y\otimes \mathrm{d}e$,\quad $\nabla 
\mathrm{d}y=\mathrm{d}x\otimes \mathrm{d}y+\mathrm{d}y\otimes \mathrm{d}x+
\mathrm{d}y\otimes \mathrm{d}y,$
\[ R_{\nabla }\mathrm{d}e=\mathrm{d}y\wedge \left( \mathrm{d}e\otimes \mathrm{d
}e\right) +\mathrm{d}e\wedge \left( \mathrm{d}x\otimes \mathrm{d}x\right)
,\quad R_{\nabla }\mathrm{d}x=\mathrm{d}x\wedge \left( \mathrm{d}e\otimes 
\mathrm{d}y+\mathrm{d}y\otimes \mathrm{d}e\right),\]
\[ R_{\nabla }\mathrm{d}y=
\mathrm{d}e\wedge \left( \mathrm{d}y\otimes \mathrm{d}y\right) +\mathrm{d}
y\wedge \left( \mathrm{d}x\otimes \mathrm{d}x\right);\]

(F.I.3) $\nabla \mathrm{d}e=\mathrm{d}e\otimes \mathrm{d}e+\mathrm{d}%
x\otimes \mathrm{d}y+\mathrm{d}y\otimes \mathrm{d}x+\mathrm{d}y\otimes 
\mathrm{d}y$,\quad $\nabla \mathrm{d}x=\mathrm{d}e\otimes \mathrm{d}e+%
\mathrm{d}y\otimes \mathrm{d}y$,\quad $\nabla \mathrm{d}y=\mathrm{d}e\otimes 
\mathrm{d}x+\mathrm{d}e\otimes \mathrm{d}y+\mathrm{d}x\otimes \mathrm{d}e+%
\mathrm{d}y\otimes \mathrm{d}e,$
\[ R_{\nabla }\mathrm{d}e=\mathrm{d}y\wedge \left( \mathrm{d}e\otimes \mathrm{d
}e\right) +\mathrm{d}e\wedge \left( \mathrm{d}x\otimes \mathrm{d}x\right)
,\quad R_{\nabla }\mathrm{d}x=\mathrm{d}x\wedge \left( \mathrm{d}e\otimes 
\mathrm{d}y+\mathrm{d}y\otimes \mathrm{d}e\right),\]
\[ R_{\nabla }\mathrm{d}y=
\mathrm{d}e\wedge \left( \mathrm{d}y\otimes \mathrm{d}y\right) +\mathrm{d}
y\wedge \left( \mathrm{d}x\otimes \mathrm{d}x\right).\]

- for metric $g_{F.II}$

(F.II.1) $\nabla \mathrm{d} e=0 $,\quad $\nabla \mathrm{d} x=\mathrm{d}
e\otimes \mathrm{d} x + \mathrm{d} e\otimes \mathrm{d} y + \mathrm{d}
x\otimes \mathrm{d} e + \mathrm{d} x\otimes \mathrm{d} y + \mathrm{d}
y\otimes \mathrm{d} e + \mathrm{d} y\otimes \mathrm{d} x $,\quad $\nabla 
\mathrm{d} y= \mathrm{d} e\otimes \mathrm{d} e + \mathrm{d} e\otimes \mathrm{
d} y + \mathrm{d} y\otimes \mathrm{d} e + \mathrm{d} y\otimes \mathrm{d} y,\quad R_\nabla=0$;

(F.II.2) $\nabla \mathrm{d} e=0 $,\quad $\nabla \mathrm{d} x=\mathrm{d}
e\otimes \mathrm{d} x + \mathrm{d} e\otimes \mathrm{d} y + \mathrm{d}
x\otimes \mathrm{d} e + \mathrm{d} x\otimes \mathrm{d} x + \mathrm{d}
y\otimes \mathrm{d} e + \mathrm{d} y\otimes \mathrm{d} y $,\quad $\nabla 
\mathrm{d} y= \mathrm{d} e\otimes \mathrm{d} x + \mathrm{d} e\otimes \mathrm{
d} y + \mathrm{d} x\otimes \mathrm{d} e + \mathrm{d} x\otimes \mathrm{d} y + 
\mathrm{d} y\otimes \mathrm{d} e + \mathrm{d} y\otimes \mathrm{d} x,\quad R_\nabla=0$;

(F.II.3) $\nabla \mathrm{d}e=0$,\quad $\nabla \mathrm{d}x=+\mathrm{d}
e\otimes \mathrm{d}e+\mathrm{d}e\otimes \mathrm{d}x+\mathrm{d}x\otimes 
\mathrm{d}e+\mathrm{d}x\otimes \mathrm{d}x$,\newline
$\nabla \mathrm{d}y=+\mathrm{d}e\otimes \mathrm{d}e+\mathrm{d}e\otimes 
\mathrm{d}y+\mathrm{d}x\otimes \mathrm{d}x+\mathrm{d}x\otimes \mathrm{d}y+
\mathrm{d}y\otimes \mathrm{d}e+\mathrm{d}y\otimes \mathrm{d}x,\quad R_\nabla=0$.

- for metric $g_{F.III}$

(F.III.1) $\nabla \mathrm{d}e=0$,\quad $\nabla \mathrm{d}x=\mathrm{d}
x\otimes \mathrm{d}x$,\newline
$\nabla \mathrm{d}y=\mathrm{d}e\otimes \mathrm{d}x+\mathrm{d}x\otimes 
\mathrm{d}e+\mathrm{d}x\otimes \mathrm{d}x+\mathrm{d}x\otimes \mathrm{d}y+
\mathrm{d}y\otimes \mathrm{d}x,\quad R_\nabla=0;$

(F.III.2) $\nabla \mathrm{d}e=\mathrm{d}e\otimes \mathrm{d}x+\mathrm{d}
x\otimes \mathrm{d}e+\mathrm{d}x\otimes \mathrm{d}x+\mathrm{d}y\otimes 
\mathrm{d}y$,\quad $\nabla \mathrm{d}x=\mathrm{d}e\otimes \mathrm{d}x+
\mathrm{d}x\otimes \mathrm{d}e+\mathrm{d}x\otimes \mathrm{d}x+\mathrm{d}
x\otimes \mathrm{d}y+\mathrm{d}y\otimes \mathrm{d}x+\mathrm{d}y\otimes 
\mathrm{d}y$,\quad $\nabla \mathrm{d}y=\mathrm{d}e\otimes \mathrm{d}y+
\mathrm{d}x\otimes \mathrm{d}x+\mathrm{d}y\otimes \mathrm{d}e+\mathrm{d}
y\otimes \mathrm{d}y,$
\[ R_{\nabla }\mathrm{d}e=\mathrm{d}x\wedge \left( \mathrm{d}e\otimes \mathrm{d
}e\right) +\mathrm{d}e\wedge \left( \mathrm{d}x\otimes \mathrm{d}x+\mathrm{d}
x\otimes \mathrm{d}y+\mathrm{d}y\otimes \mathrm{d}x\right) ,\quad \]
\[ R_{\nabla }\mathrm{d}x=\mathrm{d}x\wedge \left( \mathrm{d}e\otimes \mathrm{d
}e\right) +\mathrm{d}e\wedge \left( \mathrm{d}x\otimes \mathrm{d}x\right) +
\mathrm{d}y\wedge \left( \mathrm{d}x\otimes \mathrm{d}x\right) ,\quad\]
\[ R_{\nabla }\mathrm{d}y=\mathrm{d}y\wedge \left( \mathrm{d}e\otimes \mathrm{d
}e+\mathrm{d}e\otimes \mathrm{d}x+\mathrm{d}x\otimes \mathrm{d}e+\mathrm{d}
x\otimes \mathrm{d}x\right) +\mathrm{d}x\wedge \left( \mathrm{d}y\otimes 
\mathrm{d}y\right);\]

(F.III.3) $\nabla \mathrm{d}e=\mathrm{d}e\otimes \mathrm{d}x+\mathrm{d}
e\otimes \mathrm{d}y+\mathrm{d}x\otimes \mathrm{d}e+\mathrm{d}y\otimes 
\mathrm{d}e+\mathrm{d}y\otimes \mathrm{d}y$,\quad $\nabla \mathrm{d}x=
\mathrm{d}e\otimes \mathrm{d}x+\mathrm{d}e\otimes \mathrm{d}y+\mathrm{d}
x\otimes \mathrm{d}e+\mathrm{d}x\otimes \mathrm{d}y+\mathrm{d}y\otimes 
\mathrm{d}e+\mathrm{d}y\otimes \mathrm{d}x$,\quad $\nabla \mathrm{d}y=
\mathrm{d}e\otimes \mathrm{d}y+\mathrm{d}x\otimes \mathrm{d}y+\mathrm{d}
y\otimes \mathrm{d}e+\mathrm{d}y\otimes \mathrm{d}x,$
\[ R_{\nabla }\mathrm{d}e=\mathrm{d}x\wedge \left( \mathrm{d}e\otimes \mathrm{d
}e\right) +\mathrm{d}e\wedge \left( \mathrm{d}x\otimes \mathrm{d}x+\mathrm{d}
x\otimes \mathrm{d}y+\mathrm{d}y\otimes \mathrm{d}x\right) ,\quad \]
\[R_{\nabla }\mathrm{d}x=\mathrm{d}x\wedge \left( \mathrm{d}e\otimes \mathrm{d
}e\right) +\mathrm{d}e\wedge \left( \mathrm{d}x\otimes \mathrm{d}x\right) +
\mathrm{d}y\wedge \left( \mathrm{d}x\otimes \mathrm{d}x\right) ,\quad\]
\[R_{\nabla }\mathrm{d}y=\mathrm{d}y\wedge \left( \mathrm{d}e\otimes \mathrm{d
}e+\mathrm{d}e\otimes \mathrm{d}x+\mathrm{d}x\otimes \mathrm{d}e+\mathrm{d}
x\otimes \mathrm{d}x\right) +\mathrm{d}x\wedge \left( \mathrm{d}y\otimes 
\mathrm{d}y\right).\]

- for metric $g_{F.IV}$

(F.IV.1) $\nabla \mathrm{d}e=0$,\quad $\nabla \mathrm{d}x=\mathrm{d}e\otimes 
\mathrm{d}y+\mathrm{d}y\otimes \mathrm{d}e$,\quad $\nabla \mathrm{d}y=
\mathrm{d}e\otimes \mathrm{d}e,\quad R_\nabla=0$;

(F.IV.2) $\nabla \mathrm{d}e=\mathrm{d}e\otimes \mathrm{d}e+\mathrm{d}
e\otimes \mathrm{d}y+\mathrm{d}x\otimes \mathrm{d}x+\mathrm{d}x\otimes 
\mathrm{d}y+\mathrm{d}y\otimes \mathrm{d}e+\mathrm{d}y\otimes \mathrm{d}x+
\mathrm{d}y\otimes \mathrm{d}y$,\quad $\nabla \mathrm{d}x=\mathrm{d}x\otimes 
\mathrm{d}y+\mathrm{d}y\otimes \mathrm{d}x$,\quad $\nabla \mathrm{d}y=
\mathrm{d}e\otimes \mathrm{d}x+\mathrm{d}e\otimes \mathrm{d}y+\mathrm{d}
x\otimes \mathrm{d}e+\mathrm{d}y\otimes \mathrm{d}e,$
\[ R_{\nabla }\mathrm{d}e=\mathrm{d}x\wedge \left( \mathrm{d}e\otimes \mathrm{d
}e\right) +\mathrm{d}y\wedge \left( \mathrm{d}e\otimes \mathrm{d}e\right) +
\mathrm{d}e\wedge \left( \mathrm{d}y\otimes \mathrm{d}y\right),\]
\[ R_{\nabla }\mathrm{d}x=\mathrm{d}e\wedge \left( \mathrm{d}x\otimes \mathrm{d
}x\right) +\mathrm{d}x\wedge \left( \mathrm{d}e\otimes \mathrm{d}y+\mathrm{d}
y\otimes \mathrm{d}e+\mathrm{d}y\otimes \mathrm{d}y\right),\]
\[ R_{\nabla }\mathrm{d}y=\mathrm{d}y\wedge \left( \mathrm{d}e\otimes \mathrm{d
}x+\mathrm{d}x\otimes \mathrm{d}e\right) +\mathrm{d}e\wedge \left( \mathrm{d}
y\otimes \mathrm{d}y\right);\]

(F.IV.3) $\nabla \mathrm{d} e= +\mathrm{d} e\otimes \mathrm{d} x + \mathrm{d}
x\otimes \mathrm{d} e + \mathrm{d} x\otimes \mathrm{d} x + \mathrm{d}
y\otimes \mathrm{d} y$,\newline
$\nabla \mathrm{d} x=\mathrm{d} e\otimes \mathrm{d} x + \mathrm{d} e\otimes 
\mathrm{d} y + \mathrm{d} x\otimes \mathrm{d} e + \mathrm{d} x\otimes 
\mathrm{d} y + \mathrm{d} y\otimes \mathrm{d} e + \mathrm{d} y\otimes 
\mathrm{d} x $,\newline
$\nabla \mathrm{d} y= \mathrm{d} e\otimes \mathrm{d} x + \mathrm{d} e\otimes 
\mathrm{d} y + \mathrm{d} x\otimes \mathrm{d} e + \mathrm{d} x\otimes 
\mathrm{d} x + \mathrm{d} x\otimes \mathrm{d} y + \mathrm{d} y\otimes 
\mathrm{d} e + \mathrm{d} y\otimes \mathrm{d} x + \mathrm{d} y\otimes 
\mathrm{d} y $,
\[ R_{\nabla }\mathrm{d}e=\mathrm{d}x\wedge \left( \mathrm{d}e\otimes \mathrm{d
}e\right) +\mathrm{d}y\wedge \left( \mathrm{d}e\otimes \mathrm{d}e\right) +
\mathrm{d}e\wedge \left( \mathrm{d}y\otimes \mathrm{d}y\right),\]
\[ R_{\nabla }\mathrm{d}x=\mathrm{d}e\wedge \left( \mathrm{d}x\otimes \mathrm{d
}x\right) +\mathrm{d}x\wedge \left( \mathrm{d}e\otimes \mathrm{d}y+\mathrm{d}
y\otimes \mathrm{d}e+\mathrm{d}y\otimes \mathrm{d}y\right),\]
\[ R_{\nabla }\mathrm{d}y=\mathrm{d}y\wedge \left( \mathrm{d}e\otimes \mathrm{d
}x+\mathrm{d}x\otimes \mathrm{d}e\right) +\mathrm{d}e\wedge \left( \mathrm{d}
y\otimes \mathrm{d}y\right).\]

- for metric $g_{F.V}$

(F.V.1) $\nabla \mathrm{d}e=\mathrm{d}e\otimes \mathrm{d}y+\mathrm{d}
x\otimes \mathrm{d}x+\mathrm{d}y\otimes \mathrm{d}e$,\quad $\nabla \mathrm{d}
x=\mathrm{d}x\otimes \mathrm{d}y+\mathrm{d}y\otimes \mathrm{d}x$,\quad \newline
$\nabla \mathrm{d}y=\mathrm{d}e\otimes \mathrm{d}x+\mathrm{d}x\otimes \mathrm{
d}e+\mathrm{d}x\otimes \mathrm{d}x,$
\[ R_{\nabla }\mathrm{d}e=\mathrm{d}x\wedge \left( \mathrm{d}e\otimes \mathrm{d
}e\right) +\mathrm{d}e\wedge \left( \mathrm{d}x\otimes \mathrm{d}x+\mathrm{d}
y\otimes \mathrm{d}y\right),\]
\[ R_{\nabla }\mathrm{d}x=\mathrm{d}e\wedge \left( \mathrm{d}x\otimes \mathrm{d
}x\right) +\mathrm{d}x\wedge \left( \mathrm{d}y\otimes \mathrm{d}y\right),\]
\[ R_{\nabla }\mathrm{d}y=\mathrm{d}y\wedge \left( \mathrm{d}e\otimes \mathrm{d
}x+\mathrm{d}x\otimes \mathrm{d}e\right) +\mathrm{d}y\wedge \left( \mathrm{d}
x\otimes \mathrm{d}x\right);\]

(F.V.2) $\nabla \mathrm{d}e=\mathrm{d}e\otimes \mathrm{d}e+\mathrm{d}
x\otimes \mathrm{d}y+\mathrm{d}y\otimes \mathrm{d}x+\mathrm{d}y\otimes 
\mathrm{d}y$,\quad $\nabla \mathrm{d}x=\mathrm{d}e\otimes \mathrm{d}e+
\mathrm{d}e\otimes \mathrm{d}y+\mathrm{d}x\otimes \mathrm{d}y+\mathrm{d}
y\otimes \mathrm{d}e+\mathrm{d}y\otimes \mathrm{d}x+\mathrm{d}y\otimes 
\mathrm{d}y$,\quad $\nabla \mathrm{d}y=\mathrm{d}e\otimes \mathrm{d}e+
\mathrm{d}e\otimes \mathrm{d}x+\mathrm{d}e\otimes \mathrm{d}y+\mathrm{d}
x\otimes \mathrm{d}e+\mathrm{d}y\otimes \mathrm{d}e+\mathrm{d}y\otimes 
\mathrm{d}y,$
\[R_{\nabla }\mathrm{d}e=\mathrm{d}x\wedge \left( \mathrm{d}e\otimes \mathrm{d
}e\right) +\mathrm{d}e\wedge \left( \mathrm{d}x\otimes \mathrm{d}x+\mathrm{d}
y\otimes \mathrm{d}y\right),\]
\[ R_{\nabla }\mathrm{d}x=\mathrm{d}e\wedge \left( \mathrm{d}x\otimes \mathrm{d
}x\right) +\mathrm{d}x\wedge \left( \mathrm{d}y\otimes \mathrm{d}y\right), \]
\[ R_{\nabla }\mathrm{d}y=\mathrm{d}y\wedge \left( \mathrm{d}e\otimes \mathrm{d
}x+\mathrm{d}x\otimes \mathrm{d}e\right) +\mathrm{d}y\wedge \left( \mathrm{d}
x\otimes \mathrm{d}x\right);\]

(F.V.3) $\nabla \mathrm{d}e=0$,\quad $\nabla \mathrm{d}x=\mathrm{d}e\otimes 
\mathrm{d}e+\mathrm{d}e\otimes \mathrm{d}x+\mathrm{d}x\otimes \mathrm{d}e+
\mathrm{d}x\otimes \mathrm{d}x+\mathrm{d}y\otimes \mathrm{d}y$,\newline
$\nabla \mathrm{d}y=\mathrm{d}e\otimes \mathrm{d}y+\mathrm{d}x\otimes 
\mathrm{d}y+\mathrm{d}y\otimes \mathrm{d}e+\mathrm{d}y\otimes \mathrm{d}x,\quad R_\nabla=0$.

- for metric $g_{F.VI}$

(F.VI.1) $\nabla \mathrm{d}e=0$,\quad $\nabla \mathrm{d}x=\mathrm{d}x\otimes 
\mathrm{d}y+\mathrm{d}y\otimes \mathrm{d}x$,\quad $\nabla \mathrm{d}y=
\mathrm{d}y\otimes \mathrm{d}y,\quad R_\nabla=0$;

(F.VI.2) $\nabla \mathrm{d}e=\mathrm{d}e\otimes \mathrm{d}e+\mathrm{d}
e\otimes \mathrm{d}x+\mathrm{d}x\otimes \mathrm{d}e+\mathrm{d}x\otimes 
\mathrm{d}x+\mathrm{d}x\otimes \mathrm{d}y+\mathrm{d}y\otimes \mathrm{d}x$
,\quad $\nabla \mathrm{d}x=\mathrm{d}e\otimes \mathrm{d}e+\mathrm{d}e\otimes 
\mathrm{d}y+\mathrm{d}x\otimes \mathrm{d}x+\mathrm{d}y\otimes \mathrm{d}e$
,\quad $\nabla \mathrm{d}y=\mathrm{d}e\otimes \mathrm{d}e+\mathrm{d}e\otimes 
\mathrm{d}x+\mathrm{d}e\otimes \mathrm{d}y+\mathrm{d}x\otimes \mathrm{d}e+
\mathrm{d}x\otimes \mathrm{d}x+\mathrm{d}y\otimes \mathrm{d}e,$
\[ R_{\nabla }\mathrm{d}e=\mathrm{d}y\wedge \left( \mathrm{d}x\otimes \mathrm{d
}x\right) +\mathrm{d}e\wedge \left( \mathrm{d}y\otimes \mathrm{d}x+\mathrm{d}
x\otimes \mathrm{d}y\right),\]
\[ R_{\nabla }\mathrm{d}x=\mathrm{d}e\wedge \left( \mathrm{d}y\otimes \mathrm{d
}y\right) +\mathrm{d}y\wedge \left( \mathrm{d}x\otimes \mathrm{d}x\right),\]
\[ R_{\nabla }\mathrm{d}y=\mathrm{d}e\wedge \left( \mathrm{d}y\otimes \mathrm{d
}y+\mathrm{d}x\otimes \mathrm{d}y+\mathrm{d}y\otimes \mathrm{d}x\right) +
\mathrm{d}y\wedge \left( \mathrm{d}x\otimes \mathrm{d}x\right);\]

(F.VI.3) $\nabla \mathrm{d}e=\mathrm{d}e\otimes \mathrm{d}e+\mathrm{d}
e\otimes \mathrm{d}y+\mathrm{d}x\otimes \mathrm{d}x+\mathrm{d}x\otimes 
\mathrm{d}y+\mathrm{d}y\otimes \mathrm{d}e+\mathrm{d}y\otimes \mathrm{d}x+
\mathrm{d}y\otimes \mathrm{d}y$,\quad $\nabla \mathrm{d}x=\mathrm{d}e\otimes 
\mathrm{d}e+\mathrm{d}x\otimes \mathrm{d}x+\mathrm{d}x\otimes \mathrm{d}y+
\mathrm{d}y\otimes \mathrm{d}x+\mathrm{d}y\otimes \mathrm{d}y$,\newline
$\nabla \mathrm{d}y=\mathrm{d}e\otimes \mathrm{d}x+\mathrm{d}e\otimes 
\mathrm{d}y+\mathrm{d}x\otimes \mathrm{d}e+\mathrm{d}x\otimes \mathrm{d}y+
\mathrm{d}y\otimes \mathrm{d}e+\mathrm{d}y\otimes \mathrm{d}x,$
\[ R_{\nabla }\mathrm{d}e=\mathrm{d}e\wedge \left( \mathrm{d}y\otimes \mathrm{d
}x+\mathrm{d}x\otimes \mathrm{d}y+\mathrm{d}y\otimes \mathrm{d}y\right) ,\]
\[ R_{\nabla }\mathrm{d}x=\mathrm{d}x\wedge \left( \mathrm{d}e\otimes \mathrm{d
}e+\mathrm{d}y\otimes \mathrm{d}y\right) +\mathrm{d}y\wedge \left( \mathrm{d}
x\otimes \mathrm{d}x\right) ,\]
\[ R_{\nabla }\mathrm{d}y=\mathrm{d}y\wedge \left( \mathrm{d}e\otimes \mathrm{d
}e\right) +\mathrm{d}x\wedge \left( \mathrm{d}y\otimes \mathrm{d}y\right) .\]

- for metric $g_{F.VII}$

(F.VII.1) $\nabla \mathrm{d}e=0$,\quad $\nabla \mathrm{d}x=\mathrm{d}
e\otimes \mathrm{d}e+\mathrm{d}e\otimes \mathrm{d}x+\mathrm{d}e\otimes 
\mathrm{d}y+\mathrm{d}x\otimes \mathrm{d}e+\mathrm{d}y\otimes \mathrm{d}e$
,\newline $\nabla \mathrm{d}y=\mathrm{d}e\otimes \mathrm{d}x+\mathrm{d}e\otimes 
\mathrm{d}y+\mathrm{d}x\otimes \mathrm{d}e+\mathrm{d}y\otimes \mathrm{d}e,\quad 
R_\nabla=0$;

(F.VII.2) $\nabla \mathrm{d}e=\mathrm{d}e\otimes \mathrm{d}y+\mathrm{d}
x\otimes \mathrm{d}x+\mathrm{d}y\otimes \mathrm{d}e$,\quad
$\nabla \mathrm{d}x=\mathrm{d}e\otimes \mathrm{d}x+\mathrm{d}x\otimes 
\mathrm{d}e+\mathrm{d}x\otimes \mathrm{d}x+\mathrm{d}x\otimes \mathrm{d}y+
\mathrm{d}y\otimes \mathrm{d}x$,\quad $\nabla \mathrm{d}y=\mathrm{d}e\otimes 
\mathrm{d}y+\mathrm{d}x\otimes \mathrm{d}y+\mathrm{d}y\otimes \mathrm{d}e+
\mathrm{d}y\otimes \mathrm{d}x,$\\
\[ R_{\nabla }\mathrm{d}e=\mathrm{d}y\wedge \left( \mathrm{d}e\otimes \mathrm{d
}e\right) +\mathrm{d}e\wedge \left( \mathrm{d}x\otimes \mathrm{d}x+\mathrm{d}
x\otimes \mathrm{d}y+\mathrm{d}y\otimes \mathrm{d}x+\mathrm{d}y\otimes 
\mathrm{d}y\right) ,\]
\[ R_{\nabla }\mathrm{d}x=\mathrm{d}x\wedge \left( \mathrm{d}e\otimes \mathrm{d
}e+\mathrm{d}e\otimes \mathrm{d}y+\mathrm{d}y\otimes \mathrm{d}e\right) +
\mathrm{d}y\wedge \left( \mathrm{d}x\otimes \mathrm{d}x\right) +\mathrm{d}
x\wedge \left( \mathrm{d}y\otimes \mathrm{d}y\right) ,\]
\[ R_{\nabla }\mathrm{d}y=\mathrm{d}y\wedge \left( \mathrm{d}e\otimes \mathrm{d
}e\right) +\mathrm{d}e\wedge \left( \mathrm{d}y\otimes \mathrm{d}y+\mathrm{d}
x\otimes \mathrm{d}x\right) +\mathrm{d}y\wedge \left( \mathrm{d}e\otimes 
\mathrm{d}x+\mathrm{d}x\otimes \mathrm{d}e+\mathrm{d}x\otimes \mathrm{d}
x\right);\]

(F.VII.3) $\nabla \mathrm{d}e=\mathrm{d}e\otimes \mathrm{d}e+\mathrm{d}
e\otimes \mathrm{d}x+\mathrm{d}x\otimes \mathrm{d}e+\mathrm{d}x\otimes 
\mathrm{d}x+\mathrm{d}x\otimes \mathrm{d}y+\mathrm{d}y\otimes \mathrm{d}x$
,\quad $\nabla \mathrm{d}x=\mathrm{d}e\otimes \mathrm{d}y+\mathrm{d}y\otimes 
\mathrm{d}e$,\quad $\nabla \mathrm{d}y=\mathrm{d}e\otimes \mathrm{d}x+
\mathrm{d}e\otimes \mathrm{d}y+\mathrm{d}x\otimes \mathrm{d}e+\mathrm{d}
x\otimes \mathrm{d}y+\mathrm{d}y\otimes \mathrm{d}e+\mathrm{d}y\otimes 
\mathrm{d}x$,
\[ R_{\nabla }\mathrm{d}e=\mathrm{d}y\wedge \left( \mathrm{d}e\otimes \mathrm{d
}e\right) +\mathrm{d}e\wedge \left( \mathrm{d}x\otimes \mathrm{d}x+\mathrm{d}
x\otimes \mathrm{d}y+\mathrm{d}y\otimes \mathrm{d}x+\mathrm{d}y\otimes 
\mathrm{d}y\right),\]
\[ R_{\nabla }\mathrm{d}x=\mathrm{d}x\wedge \left( \mathrm{d}e\otimes \mathrm{d
}e+\mathrm{d}e\otimes \mathrm{d}y+\mathrm{d}y\otimes \mathrm{d}e\right) +
\mathrm{d}y\wedge \left( \mathrm{d}x\otimes \mathrm{d}x\right) +\mathrm{d}
x\wedge \left( \mathrm{d}y\otimes \mathrm{d}y\right),\]
\[ R_{\nabla }\mathrm{d}y=\mathrm{d}y\wedge \left( \mathrm{d}e\otimes \mathrm{d
}e+\mathrm{d}x\otimes \mathrm{d}x\right) +\mathrm{d}e\wedge \left( \mathrm{d}
y\otimes \mathrm{d}y\right) +\mathrm{d}x\wedge \left( \mathrm{d}y\otimes 
\mathrm{d}y\right).\]

%
%
%

\section{Partial results for $n=4$}\label{secn4}

For $n=4$ the analysis is rather more complicated  but the outlined methods work for the classification of 4-dimensional unital 
algebras and we find 16 up to isomorphism, and hence this many inner differential structures for polynomials in 4 variables over $\F_2$.
These are summarised in Table 5. 

\begin{figure}\begin{tabular}{|l|l|}
\hline
& $
\begin{array}{c}
\left[ \mathrm{d}e,e\right] =\mathrm{d}e,\quad \left[ \mathrm{d}e,x\right] =
\mathrm{d}x=\left[ \mathrm{d}x,e\right] ,\quad \left[ \mathrm{d}e,y\right] =
\mathrm{d}y=\left[ \mathrm{d}y,e\right] ,\quad \left[ \mathrm{d}e,z\right] =
\mathrm{d}z=\left[ \mathrm{d}z,e\right] 
\end{array}
$ \\ \hline
A & $
\begin{array}{c}
\left[ \mathrm{d}x,x\right] =0,\quad \left[ \mathrm{d}x,y\right] =0=\left[ 
\mathrm{d}y,x\right] ,\quad \left[ \mathrm{d}x,z\right] =0=\left[ \mathrm{d}
z,x\right]  \\ 
\left[ \mathrm{d}y,z\right] =0=\left[ \mathrm{d}z,y\right] ,\quad \left[ 
\mathrm{d}y,y\right] =0=\left[ \mathrm{d}z,z\right] 
\end{array}
$ \\ \hline
B & $
\begin{array}{c}
\left[ \mathrm{d}x,x\right] =\mathrm{d}z,\quad \left[ \mathrm{d}x,y\right]
=0=\left[ \mathrm{d}y,x\right] ,\quad \left[ \mathrm{d}x,z\right] =0=\left[ 
\mathrm{d}z,x\right]  \\ 
\left[ \mathrm{d}y,z\right] =0=\left[ \mathrm{d}z,y\right] ,\quad \left[ 
\mathrm{d}y,y\right] =0=\left[ \mathrm{d}z,z\right] 
\end{array}
$ \\ \hline
C & $~
\begin{array}{c}
\left[ \mathrm{d}x,x\right] =\mathrm{d}x,\quad \left[ \mathrm{d}x,y\right]
=0=\left[ \mathrm{d}y,x\right] ,\quad \left[ \mathrm{d}x,z\right] =0=\left[ 
\mathrm{d}z,x\right]  \\ 
\left[ \mathrm{d}y,z\right] =0=\left[ \mathrm{d}z,y\right] ,\quad \left[ 
\mathrm{d}y,y\right] =0=\left[ \mathrm{d}z,z\right] 
\end{array}
$ \\ \hline
D & $
\begin{array}{c}
\left[ \mathrm{d}x,x\right] =\mathrm{d}x,\quad \left[ \mathrm{d}x,y\right] =
\mathrm{d}y=\left[ \mathrm{d}y,x\right] ,\quad \left[ \mathrm{d}x,z\right]
=0=\left[ \mathrm{d}z,x\right]  \\ 
\left[ \mathrm{d}y,z\right] =0=\left[ \mathrm{d}z,y\right] ,\quad \left[ 
\mathrm{d}y,y\right] =0=\left[ \mathrm{d}z,z\right] 
\end{array}
$ \\ \hline
E & $~
\begin{array}{c}
\left[ \mathrm{d}x,x\right] =0,\quad \left[ \mathrm{d}x,y\right] =\mathrm{d}
z=\left[ \mathrm{d}y,x\right] ,\quad \left[ \mathrm{d}x,z\right] =0=\left[ 
\mathrm{d}z,x\right]  \\ 
\left[ \mathrm{d}y,z\right] =0=\left[ \mathrm{d}z,y\right] ,\quad \left[ 
\mathrm{d}y,y\right] =0=\left[ \mathrm{d}z,z\right] 
\end{array}
$ \\ \hline
F & $~
\begin{array}{c}
\left[ \mathrm{d}x,x\right] =\mathrm{d}z,\quad \left[ \mathrm{d}x,y\right] =
\mathrm{d}z=\left[ \mathrm{d}y,x\right] ,\quad \left[ \mathrm{d}x,z\right]
=0=\left[ \mathrm{d}z,x\right]  \\ 
\left[ \mathrm{d}y,z\right] =0=\left[ \mathrm{d}z,y\right] ,\quad \left[ 
\mathrm{d}y,y\right] =0=\left[ \mathrm{d}z,z\right] 
\end{array}
$ \\ \hline
G & $
\begin{array}{c}
\left[ \mathrm{d}x,x\right] =\mathrm{d}y,\quad \left[ \mathrm{d}x,y\right] =
\mathrm{d}z=\left[ \mathrm{d}y,x\right] ,\quad \left[ \mathrm{d}x,z\right]
=0=\left[ \mathrm{d}z,x\right]  \\ 
\left[ \mathrm{d}y,z\right] =0=\left[ \mathrm{d}z,y\right] ,\quad \left[ 
\mathrm{d}y,y\right] =0=\left[ \mathrm{d}z,z\right] 
\end{array}
$ \\ \hline
H & $
\begin{array}{c}
\left[ \mathrm{d}x,x\right] =\mathrm{d}e+\mathrm{d}x,\quad \left[ \mathrm{d}
x,y\right] =\mathrm{d}y+\mathrm{d}z=\left[ \mathrm{d}y,x\right] ,\quad \left[
\mathrm{d}x,z\right] =\mathrm{d}y=\left[ \mathrm{d}z,x\right]  \\ 
\left[ \mathrm{d}y,z\right] =0=\left[ \mathrm{d}z,y\right] ,\quad \left[ 
\mathrm{d}y,y\right] =0=\left[ \mathrm{d}z,z\right] 
\end{array}
$ \\ \hline
I & $
\begin{array}{c}
\left[ \mathrm{d}x,x\right] =\mathrm{d}y,\quad \left[ \mathrm{d}x,y\right] =
\mathrm{d}x+\mathrm{d}y=\left[ \mathrm{d}y,x\right] ,\quad \left[ \mathrm{d}
x,z\right] =0=\left[ \mathrm{d}z,x\right] ,\quad  \\ 
\left[ \mathrm{d}y,z\right] =0=\left[ \mathrm{d}z,y\right] ,\quad \left[ 
\mathrm{d}y,y\right] =\mathrm{d}x,\quad \left[ \mathrm{d}z,z\right] =0
\end{array}
$ \\ \hline
J & $
\begin{array}{c}
\left[ \mathrm{d}x,x\right] =\mathrm{d}x+\mathrm{d}z,\quad \left[ \mathrm{d}
x,y\right] =\mathrm{d}x+\mathrm{d}z=\left[ \mathrm{d}y,x\right] ,\quad \left[
\mathrm{d}x,z\right] =0=\left[ \mathrm{d}z,x\right] ,\quad  \\ 
\left[ \mathrm{d}y,z\right] =0=\left[ \mathrm{d}z,y\right] ,\quad \left[ 
\mathrm{d}y,y\right] =\mathrm{d}x,\quad \left[ \mathrm{d}z,z\right] =0
\end{array}
$ \\ \hline
K & $
\begin{array}{c}
\left[ \mathrm{d}x,x\right] =\mathrm{d}x,\quad \left[ \mathrm{d}x,y\right]
=0=\left[ \mathrm{d}y,x\right] ,\quad \left[ \mathrm{d}x,z\right] =0=\left[ 
\mathrm{d}z,x\right]  \\ 
\left[ \mathrm{d}y,z\right] =0=\left[ \mathrm{d}z,y\right] ,\quad \left[ 
\mathrm{d}y,y\right] =\mathrm{d}y,\quad \left[ \mathrm{d}z,z\right] =0
\end{array}
$ \\ \hline
L & $
\begin{array}{c}
\left[ \mathrm{d}x,x\right] =\mathrm{d}z,\quad \left[ \mathrm{d}x,y\right]
=0=\left[ \mathrm{d}y,x\right] ,\quad \left[ \mathrm{d}x,z\right] =\mathrm{d}
e+\mathrm{d}y=\left[ \mathrm{d}z,x\right]  \\ 
\left[ \mathrm{d}y,y\right] =\mathrm{d}y,\quad \left[ \mathrm{d}y,z\right]
=0=\left[ \mathrm{d}z,y\right] ,\quad \left[ \mathrm{d}z,z\right] =\mathrm{d}
x
\end{array}
$ \\ \hline
M & $
\begin{array}{c}
\left[ \mathrm{d}x,x\right] =\mathrm{d}e+\mathrm{d}x+\mathrm{d}y+\mathrm{d}
z,\quad \left[ \mathrm{d}x,y\right] =0=\left[ \mathrm{d}y,x\right] ,\quad 
\left[ \mathrm{d}x,z\right] =\mathrm{d}e+\mathrm{d}x+\mathrm{d}y=\left[ 
\mathrm{d}z,x\right]  \\ 
\left[ \mathrm{d}y,z\right] =0=\left[ \mathrm{d}z,y\right] ,\quad \left[ 
\mathrm{d}y,y\right] =\mathrm{d}y,\quad \left[ \mathrm{d}z,z\right] =\mathrm{
d}x
\end{array}
$ \\ \hline
N & $
\begin{array}{c}
\left[ \mathrm{d}x,x\right] =\mathrm{d}z,\quad \left[ \mathrm{d}x,y\right]
=0=\left[ \mathrm{d}y,x\right] ,\quad \left[ \mathrm{d}y,z\right] =0=\left[ 
\mathrm{d}z,y\right]  \\ 
\left[ \mathrm{d}x,z\right] =\mathrm{d}x+\mathrm{d}z=\left[ \mathrm{d}z,x
\right] ,\quad \left[ \mathrm{d}y,y\right] =\mathrm{d}e+\mathrm{d}x+\mathrm{d
}y+\mathrm{d}z,\quad \left[ \mathrm{d}z,z\right] =\mathrm{d}x
\end{array}
$ \\ \hline
O & $
\begin{array}{c}
\left[ \mathrm{d}x,x\right] =\mathrm{d}e+\mathrm{d}z,\quad \left[ \mathrm{d}
x,y\right] =\mathrm{d}z=\left[ \mathrm{d}y,x\right] ,\quad \left[ \mathrm{d}
x,z\right] =\mathrm{d}e+\mathrm{d}y=\left[ \mathrm{d}z,x\right]  \\ 
\left[ \mathrm{d}y,z\right] =\mathrm{d}e=\left[ \mathrm{d}z,y\right] ,\quad 
\left[ \mathrm{d}y,y\right] =\mathrm{d}x+\mathrm{d}y,\quad \left[ \mathrm{d}
z,z\right] =\mathrm{d}e
\end{array}
$ \\ \hline
P & $
\begin{array}{c}
\left[ \mathrm{d}x,x\right] =\mathrm{d}x,\quad \left[ \mathrm{d}x,y\right]
=0=\left[ \mathrm{d}y,x\right] ,\quad \left[ \mathrm{d}x,z\right] =0=\left[ 
\mathrm{d}z,x\right]  \\ 
\left[ \mathrm{d}y,z\right] =0=\left[ \mathrm{d}z,y\right] ,\quad \left[ 
\mathrm{d}y,y\right] =\mathrm{d}y,\quad \left[ \mathrm{d}z,z\right] =\mathrm{
d}z
\end{array}
$ \\ \hline
\end{tabular}
Table 5. All possible inner differential structures for on $\F_2[e,x,y,z].$ \end{figure}

The methods are the same we have seen for $n=2,3$ so we will be brief. 
For the inner case with $\theta=\extd x^1$ by computer we get 5216 solutions to eqs. \eqref{prelie1}-\eqref{prelie2} for $\F_2[x^1,x^2,x^3,x^4]$. The isomorphisms group $G$ has the order $|G|=20160$. Checking the isomorphisms between all of the solutions of the inner case with $\theta=dx^1$ there are only 16 inequivalent differential calculi. The remaining possibilities for $\theta$ are isomorphic to the one with $\theta=\extd x^1$. We renamed $\theta=\extd e$ and the remaining variables as $x^2=x,\ x^3=y,\ x^4=z$ and listed the calculi in the table

Most of these calculi have metrics leading to too many geometries to study so explicitly as we did for $n=2,3$, so we focus on some that fit with the general examples (i)-(iii) in Section~\ref{secformalism}. 

\subsection*{(i) $V={\mathbb{F}}(X)$} We carefully make the same change of variable noted for $n=2$ and $n=3$ above (case B in those tables but case P in the $n=4$ table) to determine the quantum Levi-Civita connection in the general $x^\mu$ coordinate system where our basis of delta-functions on $X$ is labelled by $\mu=1,2,\cdots,n=|X|$ (or abstractly by $\mu\in X$ as indexing set). The change of coordinate only concerns $e=\sum_\mu x^\mu$ as the algebra identity element. For $n=2$ we let $t=e+x$ then $x^\mu=t,x$ are like spacetime coordinates. Then the metric  and quantum Levi-Civita connection found in Section~\ref{secn2} are
\[ g_B=\extd t\tens\extd t+\extd x\tens \extd x,\quad \nabla\extd t=\nabla\extd x=(\extd t+\extd x)\tens(\extd t+\extd x).\]
Similarly for $n=3$ we let $t=e+x+y$ then $x^\mu=t,x,y$ are like spacetime coordinates and the metric and three quantum Levi-Civita connection found in Section~\ref{secn3} are  $g_B=\extd t\tens\extd t+\extd x\tens\extd x+\extd y\tens\extd y$ as expected and: 

(B.1) $\nabla\extd x=0$, $\nabla \extd t=\nabla\extd x=(\extd t+\extd x)\tens (\extd t+\extd x)$;

(B.2) $\nabla\extd y=0$, $\nabla \extd t=\nabla\extd y=(\extd t+\extd y)\tens (\extd t+\extd y)$;

(B.3) $\nabla\extd t=0$, $\nabla \extd x=\nabla\extd y=(\extd x+\extd y)\tens (\extd x+\extd y)$.

From these formulae we can now extrapolate from $n=2,3$ to a general construction for quantum Levi-Civita connections for this case.

\begin{proposition}\label{Xqlc} For calculus defined by  $V=\F_2[X]$ we partition $X$ into the form $X=T\sqcup S\sqcup \bar S$ and fix a bijection $\overline{\phantom{s}}:S\to\bar S$ and letter the indices accordingly. Every such partition and bijection datum defines a quanutm Levi-Civita bimodule connection for the Euclidean metric 
\[ g=\sum_{t\in T}\extd x^t\tens \extd x^t+\sum_{s\in S}(\extd x^s\tens\extd x^s+\extd x^{\bar s}\tens\extd x^{\bar s})\]
 namely
\[ \nabla  \extd x^t=0,\quad \forall t\in T,\quad \nabla\extd x^s=\nabla\extd x^{\bar s}=(\extd x^s+\extd x^{\bar s})\tens (\extd x^s+\extd x^{\bar s}),\quad \forall s\in S\]
\[ \sigma(\extd x^s\tens\extd x^s)=\extd x^{\bar s}\tens\extd x^{\bar s},\quad \sigma(\extd x^{\bar s}\tens\extd x^{\bar s})=\extd x^s\tens\extd x^s\]
\[ \sigma(\extd x^s\tens\extd x^{\bar s})=\extd x^s\tens\extd x^{\bar s},\quad \sigma(\extd x^{\bar s}\tens\extd x^s)=\extd x^{\bar s}\tens\extd x^s\]
and otherwise the flip map. These connections have zero curvature.
\end{proposition}
\proof Once we obtained the formula based on our computer results for $n=2,3$ it is not hard to verify directly that this is quantum metric compatible and quantum torsion free from the general form of the commutation relations $[\extd x^\mu,\extd x^\nu]=\delta_{\mu\nu}\extd x^\mu$. The $\sigma$ is then uniquely determined from $\nabla$ and comes out as stated. In general we need $m=|T|$ the same cardinality as $|X|$ and there are then $[{n\atop m}]$ choices for which indices to put in $T$, followed by $(n-m-1)!!$ choices for how to pair off the remaining elements into two groups. Note that if we do not care about the labelling of indices (geometrically these are all equivalent) then we have just the integer part of $n/2$ choices for $m$ for the number of nonzero connections, but in our tables we have been distinguishing these. The zero curvature is immediate from the formulae for $\nabla$. 
 \endproof

 For $n=2$ we can take $m=2$ (the zero connection) or $m=0$ with one choice for the connection in this case. For $n=3$ we can take $m=3$ (the zero connection) or $m=1$ which three choices for which element to take for $T$ and then a unique connection for each choice. This agrees with the results from the previous tables. For $n=4$ of interest here (case P in Table~5) we can take $m=4$ (the zero connection), $m=2$ which has 6 choice for $T$ and for each of these a unique connection, or $m=0$ with $3$ choices for how to pair off the 4 indices.  For an example of the latter we can  take 
 \begin{equation}\label{n4ex} \nabla\extd t=\nabla\extd x=(\extd t+\extd x)\tens (\extd t+\extd x),\quad \nabla\extd y=\nabla\extd z=(\extd y+\extd z)\tens (\extd y+\extd z)\end{equation}
 if we write $x^\mu=t,x,y,z$.

\subsection*{(ii) $V=\F_2{\mathbb{Z}}_n$} For $n=4$ this is case G in Table~5 after a change of variables to $x^0=e, x^1=e+x, x^2=e+y, x^3=e+x+y+x$ (on the left are labels not exponents, albeit exponents with the $\circ$ product). The same methods as above for $n=2,3$ give us 8 quantum metrics for $n=4$ with matrices, written in
basis order $\extd x^1,\extd x^2,\extd x^3,\extd x^0$,
\begin{equation*}
\begin{bmatrix}
0 & 0 & 0 & 1 \\ 
0 & 0 & 1 & 0 \\ 
0 & 1 & 0 & 0 \\ 
1 & 0 & 0 & 0
\end{bmatrix}
,\quad 
\begin{bmatrix}
0 & 0 & 1 & 0 \\ 
0 & 1 & 0 & 0 \\ 
1 & 0 & 0 & 0 \\ 
0 & 0 & 0 & 1
\end{bmatrix}
,\quad 
\begin{bmatrix}
0 & 1 & 0 & 0 \\ 
1 & 0 & 0 & 0 \\ 
0 & 0 & 0 & 1 \\ 
0 & 0 & 1 & 0
\end{bmatrix}
,\quad 
\begin{bmatrix}
1 & 0 & 0 & 0 \\ 
0 & 0 & 0 & 1 \\ 
0 & 0 & 1 & 0 \\ 
0 & 1 & 0 & 0
\end{bmatrix}
,\quad 
\begin{bmatrix}
1 & 1 & 1 & 0 \\ 
1 & 1 & 0 & 1 \\ 
1 & 0 & 1 & 1 \\ 
0 & 1 & 1 & 1
\end{bmatrix}
,\quad 
\begin{bmatrix}
1 & 1 & 0 & 1 \\ 
1 & 0 & 1 & 1 \\ 
0 & 1 & 1 & 1 \\ 
1 & 1 & 1 & 0
\end{bmatrix}
,\quad 
\begin{bmatrix}
1 & 0 & 1 & 1 \\ 
0 & 1 & 1 & 1 \\ 
1 & 1 & 1 & 0 \\ 
1 & 1 & 0 & 1
\end{bmatrix}
,\quad 
\begin{bmatrix}
0 & 1 & 1 & 1 \\ 
1 & 1 & 1 & 0 \\ 
1 & 1 & 0 & 1 \\ 
1 & 0 & 1 & 1
\end{bmatrix}
\end{equation*}
of which the first four are in the general example (ii) in Section~\ref{secformalism} (for $m=1,0,3,2$ respectively). The other four are their complementary metrics with de-Morgan dual coefficients which for $n=4$ are distinct and nondegenerate. Thus the general construction together with duals gives all metrics at least for $n\le 4$ (and we suppose for all $n$). Experience with $n=2,3$ tells us to expect more than one nonzero quantum Levi-Civita connection for each metric when $n=4$.

\subsection*{(iii)  $V=A_2$} This appears as case L in Table~ 5 after a change of variables to $x^0=e, x^1=x, x^2=z$ and $x^3=e+y$ (on the left are labels not exponents, albeit exponents for the $\circ$ product) to match the basis in the general example~(iii) in Section~\ref{secformalism}. In this basis the relations for the calculus on $\F_2[x^0,x^1,x^2,x^3]$ are
\begin{equation*}
\lbrack \mathrm{d}x^0,x^0]=x^0,\quad \lbrack \mathrm{d}x^0,x^{i}]=[\mathrm{d}
x^{i},x^0]=\mathrm{d}x^{i},\quad \lbrack \mathrm{d}x^{1},x^{1}]=[\mathrm{d}
x^{2},x^{3}]=[\mathrm{d}x^{3},x^{2}]=\mathrm{d}x^{2}
\end{equation*}
\begin{equation*}
\lbrack \mathrm{d}x^{2},x^{1}]=[\mathrm{d}x^{1},x^{2}]=[\mathrm{d}
x^{3},x^{3}]=\mathrm{d}x^{3},\quad \lbrack \mathrm{d}x^{3},x^{1}]=[\mathrm{d}
x^{2},x^{2}]=[\mathrm{d}x^{1},x^{3}]=\mathrm{d}x^{1}
\end{equation*}
and by computer one has three quantum metrics, 

$g_I=\mathrm{d} x^1\otimes\mathrm{d} x^3+x^3\otimes x^1+\mathrm{d}
x^2\otimes \mathrm{d} x^2 + \mathrm{d} x^3\otimes \mathrm{d} x^3 + \mathrm{d}
x^3\otimes \mathrm{d} x^0 +\mathrm{d} x^0\otimes \mathrm{d} x^3 +\mathrm{d}
x^0\otimes \mathrm{d} x^0$

$g_{II}=\mathrm{d} x^1\otimes \mathrm{d} x^2 + \mathrm{d} x^2\otimes 
\mathrm{d} x^1 + \mathrm{d} x^3\otimes \mathrm{d} x^0 +\mathrm{d} x^0\otimes 
\mathrm{d} x^3 +\mathrm{d} e\otimes \mathrm{d} x^0$

$g_{III}=\mathrm{d} x^1\otimes \mathrm{d} x^1 + \mathrm{d} x^2\otimes 
\mathrm{d} x^3 + \mathrm{d} x^3\otimes \mathrm{d} x^2 + \mathrm{d}
x^3\otimes \mathrm{d} x^3 + \mathrm{d} x^3\otimes \mathrm{d} x^0 +\mathrm{d}
x^0\otimes \mathrm{d} x^3 +\mathrm{d} x^0\otimes \mathrm{d} x^0$

In both cases (ii),(iii) one can clearly go ahead and look for quantum Levi-Cita connections but 
would need a more powerful computer. 

\section{Discussion}\label{secdis}

In this paper our main focus was on inner differential calculi and as such we classified all noncommutative Riemannian geometries on $\F_2[x^1,\cdots, x^n]$ i.e. in $n$-dimensions and with constant coefficients, for $n\le 3$ and some results for $n=4$ or higher. There are several remarks to be made. 

First of all, the inner case was a useful restriction which is typical of strictly noncommutative geometries, taking the view that classical geometry is a somewhat special and unrepresentative limit. However, a similar analysis and classification can be done without this requirement, it just produces many more calculi. For example, for $n=2$ we find two additional families (to the three given in Table 1) namely 
$D: [\extd e,e]=e$ and $E: [\extd e,e]=\extd x$ for the non-zero commutors. 
For calculus D there exists one quantum metric and for calculus E there exist two quantum metrics. All quantum metrics and quantum Levi-Civita connections (parametrised by $\alpha,\beta\in\F_2$) are shown in Table 6. These are in addition to the classical, commutative, calculus which has the zero algebra (all products zero). 

\begin{figure}
\small{
\begin{tabular}{|l|l|l|}
\hline
& Relations & Quantum metrics and QLCs \\ \hline
D & $
\begin{array}{c}
\left[ \mathrm{d}e,e\right] =e\quad  \\ 
\left[ \mathrm{d}x,e\right] =0=\left[ \mathrm{d}e,x\right]  \\ 
\quad \left[ \mathrm{d}x,x\right] =0
\end{array}
$ & $
\begin{array}{c}
g_{D}=\mathrm{d}e\otimes \mathrm{d}e+\mathrm{d}x\otimes \mathrm{d}x \\ 
\nabla \mathrm{d}e=\alpha \mathrm{d}x\otimes \mathrm{d}e,\quad \nabla 
\mathrm{d}x=\beta \mathrm{d}x\otimes \mathrm{d}x \\ 
\nabla \mathrm{d}e=\alpha \mathrm{d}x\otimes \mathrm{d}e+\mathrm{d}x\otimes 
\mathrm{d}x,\quad \nabla \mathrm{d}x=\mathrm{d}x\otimes \mathrm{d}e+\beta 
\mathrm{d}x\otimes \mathrm{d}x \\ 
\nabla \mathrm{d}e=\mathrm{d}e\otimes \mathrm{d}x+\alpha \mathrm{d}x\otimes 
\mathrm{d}e,\quad \nabla \mathrm{d}x=\mathrm{d}e\otimes \mathrm{d}e+\beta 
\mathrm{d}x\otimes \mathrm{d}x
\end{array}
$ \\ \hline
E & $
\begin{array}{c}
\left[ \mathrm{d}e,e\right] =x,\quad  \\ 
\left[ \mathrm{d}x,e\right] =0=\left[ \mathrm{d}e,x\right]  \\ 
\quad \left[ \mathrm{d}x,x\right] =0
\end{array}
$ & $
\begin{array}{c}
g_{E.I}=\mathrm{d}e\otimes \mathrm{d}x+\mathrm{d}x\otimes \mathrm{d}e\newline
\\ 
\nabla \mathrm{d}e=\alpha \mathrm{d}x\otimes \mathrm{d}x,\quad \nabla 
\mathrm{d}x=\beta \mathrm{d}e\otimes \mathrm{d}e \\ 
\nabla \mathrm{d}e=\mathrm{d}e\otimes \mathrm{d}x+\alpha \mathrm{d}x\otimes 
\mathrm{d}x,\quad \nabla \mathrm{d}x=\beta \mathrm{d}e\otimes \mathrm{d}e \\ 
\nabla \mathrm{d}e=\mathrm{d}e\otimes \mathrm{d}e+\alpha \mathrm{d}x\otimes 
\mathrm{d}x,\quad \nabla \mathrm{d}x=\beta \mathrm{d}e\otimes \mathrm{d}e+
\mathrm{d}e\otimes \mathrm{d}x+\mathrm{d}x\otimes \mathrm{d}e \\ 
\nabla \mathrm{d}e=\mathrm{d}x\otimes \mathrm{d}e+\alpha \mathrm{d}x\otimes 
\mathrm{d}x,\quad \nabla \mathrm{d}x=\beta \mathrm{d}e\otimes \mathrm{d}e+
\mathrm{d}x\otimes \mathrm{d}x \\ 
\nabla \mathrm{d}e=\mathrm{d}e\otimes \mathrm{d}x+\mathrm{d}x\otimes \mathrm{
d}e+\alpha \mathrm{d}x\otimes \mathrm{d}x,\quad \nabla \mathrm{d}x=\beta 
\mathrm{d}e\otimes \mathrm{d}e+\mathrm{d}x\otimes \mathrm{d}x \\ 
\nabla \mathrm{d}e=\mathrm{d}e\otimes \mathrm{d}e+\mathrm{d}x\otimes \mathrm{
d}e+\alpha \mathrm{d}x\otimes \mathrm{d}x,\quad \nabla \mathrm{d}x=\beta 
\mathrm{d}e\otimes \mathrm{d}e++\mathrm{d}e\otimes \mathrm{d}x+\mathrm{d}
x\otimes \mathrm{d}e+\mathrm{d}x\otimes \mathrm{d}x \\ 
\hline
\begin{array}{c}
g_{E.II}=\mathrm{d}e\otimes \mathrm{d}x+\mathrm{d}x\otimes \mathrm{d}e+
\mathrm{d}x\otimes \mathrm{d}x \\ 
\nabla \mathrm{d}e=\alpha \mathrm{d}x\otimes \mathrm{d}x,\quad \nabla 
\mathrm{d}x=\beta \mathrm{d}e\otimes \mathrm{d}e \\ 
\nabla \mathrm{d}e=\mathrm{d}e\otimes \mathrm{d}x+\alpha \mathrm{d}x\otimes 
\mathrm{d}x,\quad \nabla \mathrm{d}x=\beta \mathrm{d}e\otimes \mathrm{d}e \\ 
\nabla \mathrm{d}e=\mathrm{d}e\otimes \mathrm{d}e+\mathrm{d}e\otimes \mathrm{
d}x+\alpha \mathrm{d}x\otimes \mathrm{d}x,\quad \nabla \mathrm{d}x=\beta 
\mathrm{d}e\otimes \mathrm{d}e+\mathrm{d}e\otimes \mathrm{d}x+\mathrm{d}
x\otimes \mathrm{d}e+\mathrm{d}x\otimes \mathrm{d}x \\ 
\nabla \mathrm{d}e=\mathrm{d}x\otimes \mathrm{d}e+\alpha \mathrm{d}x\otimes 
\mathrm{d}x,\quad \nabla \mathrm{d}x=\beta \mathrm{d}e\otimes \mathrm{d}e+
\mathrm{d}x\otimes \mathrm{d}x \\ 
\nabla \mathrm{d}e=\mathrm{d}e\otimes \mathrm{d}x+\mathrm{d}x\otimes \mathrm{
d}e+\alpha \mathrm{d}x\otimes \mathrm{d}x,\quad \nabla \mathrm{d}x=\beta 
\mathrm{d}e\otimes \mathrm{d}e+\mathrm{d}x\otimes \mathrm{d}x \\ 
\nabla \mathrm{d}e=\mathrm{d}e\otimes \mathrm{d}e+\mathrm{d}e\otimes \mathrm{
d}x+\mathrm{d}x\otimes \mathrm{d}e+\alpha \mathrm{d}x\otimes \mathrm{d}
x,\quad \nabla \mathrm{d}x=\beta \mathrm{d}e\otimes \mathrm{d}e+\mathrm{d}
e\otimes \mathrm{d}x+\mathrm{d}x\otimes \mathrm{d}e
\end{array}
\end{array}
$ \\ \hline
\end{tabular}
}
\normalsize
Table 6. All possible non-inner noncommutative geometries on $\F_2[e,x]$. Here $\alpha,\beta\in \F_2$ are parameters.
\end{figure}

Next, we should note that while our `coordinate algebra' $A=\F_2[x^1,\cdots, x^n]$ has been classical, the same formulae for differential geometries hold identically if we have commutation relations of Heisenberg/Clifford type (there being no difference over $\F_2$) defined by some $\Theta_{\mu\nu}$. For example in $n=4$ the structure of the connection for the calculus given by $V=\F_2[{\rm 4\ points}]$ suggests pair-wise grouping with relations
\[ xt+tx=1,\quad yz+zy=1\]
for an algebra $A$ of `quantum spacetime'.  The geometry is not affected as explained in Section~\ref{secformalism} since $\extd 1=0$ and since our formulae have constant coefficients. 

This brings us to a main limitation of the paper, namely to constant coefficients in the metric and connection in our $\extd x^\mu$ basis. This means that our geometries are in some sense `flat space' and indeed we checked that many of them have zero curvature. What is surprising is that even so there are so many rich possibilities for the quantum Levi-Civita connection other than $\nabla\extd x^\mu=0$ and $\sigma={\rm flip}$ which is the obvious `flat' connection, and indeed some of them are even curved. This non-uniqueness of the torsion free metric compatible bimodule connection for a given metric is also seen in some other noncommutative models, such as \cite{BegMa2}. It is also remarkable that we can't take any constant coefficients for the metric, which is a rigidity phenomenon for noncommutative geometry again seen in other models \cite{BegMa2,MaTao1}. In our case the number of metrics is far less than the potentially $2^{n(n+1)/2}$ possible coefficients values and gave our rich classification. 

We now consider applications of such noncommutative geometries. Our motivation here is that they model quantum spacetime, but one could also apply them in many other contexts such as `digital' models of quantum mechanics phase spaces or other `geometric' applications in  engineering. Apart from enumerating the different geometries (which would be relevant to a sector of quantum gravity where we sum over geometries) we can generally explore particles and fields on each noncommutative-geometric background, for example solutions of wave equations and Maxwell equations. Here the natural scalar Laplacian in our approach to noncommutative geometry is defined by $\square=(\ ,)\nabla\extd$ where $(\ ,\ ):\Omega^1\tens_A\Omega^1\to A$ is the inverse metric \cite{BegMa2}. To make this concrete we take the differential calculus be given by $V=\F_2(X)$ for a finite indexing set $X$ (so spacetime coordinates are $x^\mu$ where $\mu\in X$). We are then forced to the Euclidean metric and have connections as in Proposition~\ref{Xqlc}. The non-commutation relations in each variable and the fact that they mutually commute gives us
\[ \extd f(x^1,\cdots,x^n)=\sum_\mu (\del_\mu f)\extd x^\mu,\quad \del_\mu f(x^1,\cdots,x^n)=f(x^1,\cdots,x^\mu+1,\cdots,x^n)-f(x^1,\cdots,x^n)\]
i.e. the partial derivatives are finite difference operators. Then the Leibniz properties of a connection and evaluation against the inner product give
\[ \square f(x^1,\cdots,x^n)=\sum_\mu \del_\mu\del_\mu f\]
independently of the connection (this is because $(\extd x^s+\extd x^{\bar s},\extd x^s+\extd x^{\bar s})=0$ over $\F_2$). For example, we have
\[ f=\sum_{i_1,\cdots,i_n}a_{i_1i_2\cdots i_n}(x^1)^{2^{i_1}}\cdots (x^n)^{2^{i_n}},\quad \del^\mu f=\sum_{i_1,\cdots,i_n}a_{i_1i_2\cdots i_n}(x^1)^{2^{i_1}}\cdots\widehat{(x^\mu)^{2^{i_\mu}}}\cdots (x^n)^{2^{i_n}}\]
where we leave out the $x^\mu$. Hence such functions are  automatically zero modes of $\square$. This is probably the simplest example; other $V$ will lead to other commutation relations and other geometries. The properties and applications of such geometric wave operators would be an interesting topic for further work. Our idea is that such equations could be used to propagate information much as in a quantum computer but here modelled with `digital geometry'.  

A more speculative direction to be explored here is to use the above as a model of classical and quantum field theory. Thus we have been thinking of $A=\F_2[V]$ by which we mean polynomials in generators $x^\mu$ arising as a basis of a commutative algebra with vector space $V$. But what if $V$ is actually the spacetime coordinate algebra? For example if $V=\F_2(X)$ and $X$ is a discrete spacetime with basis $x^\mu=\delta_\mu$ of delta-functions at $\mu\in X$  then  $A$ would be functionals on $X$ (i.e. functions on the vector space of functions on $X$). We can also allow these functionals to have Heisenberg-type commutation relations as explained in \ref{secn2} above and if we are not interested in a metric then we do not need to work over $\F_2$ so we can get closer to conventional classical or quantum field theory. What we see, however, is that in this context the structure of $V$ (the classical spacetime geometry in some sense) determines a noncommutative differential on the algebra of functionals, i.e a noncommutative variational calculus. It would be interesting to reformulate such things as noncommutative Euler-Lagrange equations from this point of view now much more tied to the classical spacetime geometry. This applies even if we keep $A$ classical, i.e. are studying noncommutative variations or differentials of classical fields on $X$.

This brings us to a different classification problem. If we are interested in commutative algebras $V$ as `spacetime coordinate algebras' with the classical or possibly quantum field theory interpretation of $A$, then we should also be interested in the differential geometry of $V$ as our spacetime differential geometry. This will then connect through to any variational field equations just as it does classically when $V=C^\infty(M)$.
We again can make things simpler by letting $V$ be finite-dimensional (a finite geometry) and going `digital' by working over $\F_2$. We have already done part of the classification since we classified unital algebras over $\F_2$ up to dimension 4 (the latter being too numerous to list explicitly in the present paper). Beyond this we should consider noncommutative differential and Riemannian structures over each $V$, which is a classification problem we will address by computer methods similar to the above, in \cite{MaPac}. Such a `finite digital geometry' was obtained for the 4-dimensional algebra  $A_2$ in \cite[Prop. 5.7]{BasMa}, where it shown that there is a natural 2-dimensional differential calculus with 3 possible metrics with constant coefficients, and for each of these the paper found one quantum Levi-Civita connection other than the zero one, with zero curvature. 

For the three unital algebras $V$ of dimension 2 identified in Section~\ref{secn2} we have only the zero calculus or the universal calculus of the maximal dimension $n-1$, i.e. 1-dimensional over the algebra. The relations of the latter for each algebra are obtained by applying $\extd$ to the algebra relations, giving 

A:\quad    $[\extd x,x]=0$;\quad B: $[\extd x,x]=\extd x$;\quad  C:\quad   $[\extd x,x]=\extd x$

for the three $\Omega^1(V)$, along with $e=1$ central and killed by $\extd$. In each case we have $g=\extd x\tens\extd x$ as the only  metric and $\nabla\extd x=0$ or $\nabla\extd x=\extd x\tens\extd x$ as quantum Levi-Civita connections. For $n>2$ we have different differential structures form the zero up to the universal of dimension $n-1$ over the algebra with more nontrivial geometries arising. 

We can also allow our algebras to be noncommutative and look for other algebraic structures including Hopf algebras and solutions of the Yang-Baxter or braid relations over $\F_2$. The nice thing about doing algebra over $\F_2$ is that it could in theory be realised both in software machine code or indeed in actual silicon by means of logic gates, as follows. Apart from the motivation given, geometric elements of quantum computing could then be implemented digitally while possibly keeping some of the benefits. 

First, we can represent a vector space $V$  of dimension $m$ and fixed basis $\{e_i\}_{i=1}^{i=m}$  by a ribbon cable of $m$ wires. Then any $v\in V$ corresponds to a pattern of 0s and 1s in the wires according to $v=\sum v_i e_i$ where $v_i$ is the digital signal (0 or 1) in the $i$'th wire. Thus $v$ has $e_i$ each time there is a 1 in the $i$-the wire. Equivalently, the binary number $v_{1_1}\cdots v_{i_m}$ represents the vector $v\in V$ (there are $2^m$ states of each).   If $W$ similarly has basis $\{f_j\}$ where $j=1,\cdots,p$ then we identify $V\otimes W$ with $mp$ wires by basis  $E_{i,j}=e_i\tens f_j$, which we can organise as $m$ bundles of $p$-wire cables (one could imagine them stacked below each other). Direct sum $V\oplus W$ corresponds to a $m+p$-wire cable given by placing the $m$-wire cable for $V$ next to that of $W$, as does $V\times W$.  Algebraic operations can then be written as digitial gates assigning to all input truth table values an output truth table, as expressed in the state of the wires. 

\begin{figure}
\[\includegraphics[scale=0.4]{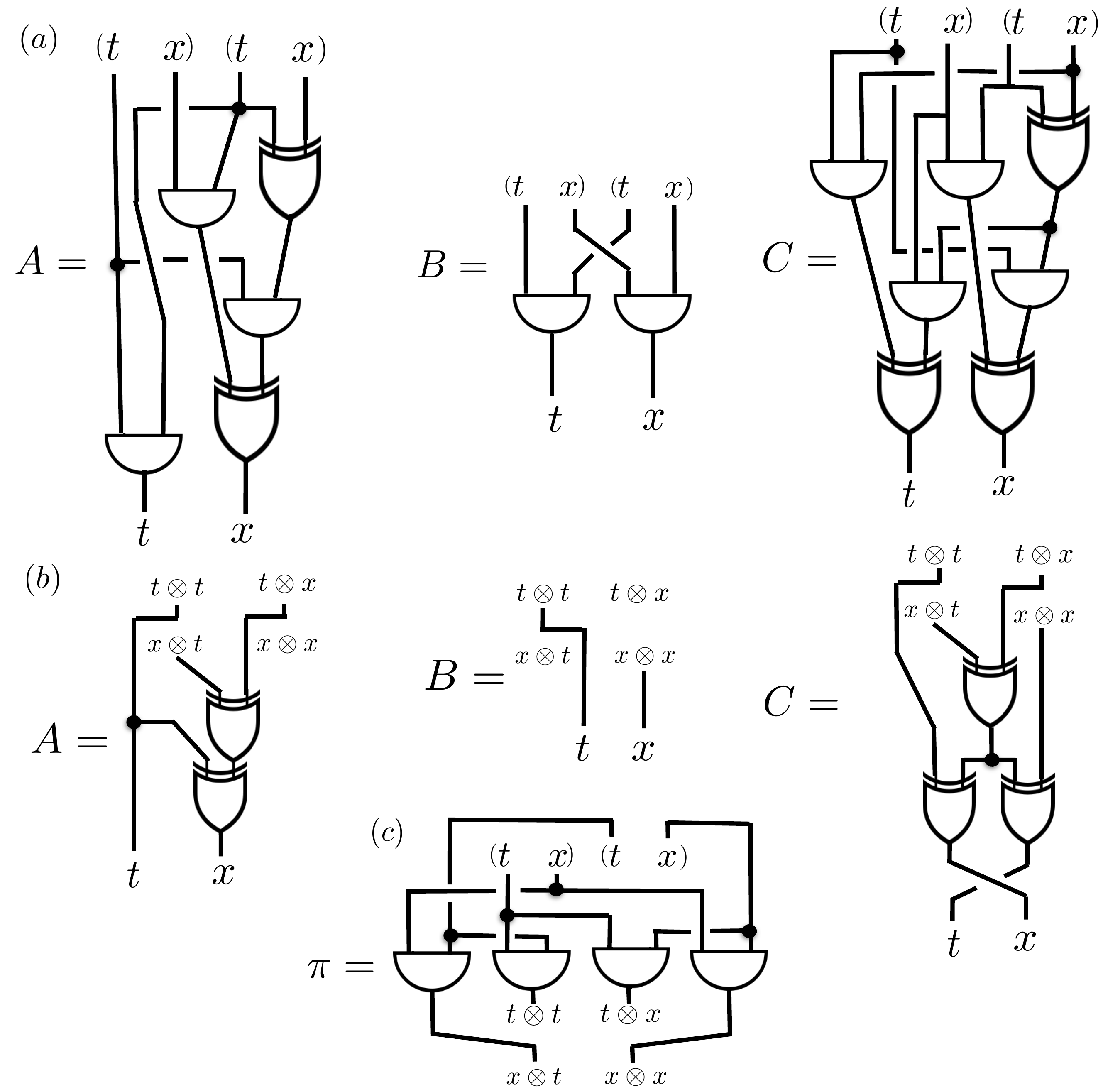}\]
\caption{Electronic circuit diagrams for each of the 3 unital products A,B,C of dimension $2$ over $\F_2$, (a) as $V\times V\to V$ and (b) as  $V\tens V\to V$. (c) shows the canonical map $V\times V\to V\tens V$}
\end{figure}

To illustrate this we show the circuit diagrams for the three $n=2$ unital algebras over $\F_2$ (cases A,B,C in Section~\ref{secn2}). Here $V$ is 2-dimensional so has 2 wires.  We chose basis $t=e+x,x$ for $V$ (this is arbitrary but recall that in the example of $\F_2[{\rm 2\ points}]$ we had $t$ and $x$ as the natural basis of $\delta$-functions for the two points). Thus the correspondence we use between vectors and digital states will be
\[ 0=00,\quad x=01,\quad t=e+x=10,\quad e=11.\]
The easiest representation is to define the algebra product as a map $V\times V\to V$ with two 2-wire inputs one for each element of $V$ and a 2-wire output which give the 3 different algebra products shown in part (a). The flat edged `product' denotes AND which is 1 exactly when both inputs are. The other curve-edged `product' denotes symmetric difference or XOR (exclusive OR) which is 1 exactly when the two inputs are different. The desired outcomes can be expressed as Boolean algebra or more precisely as a Boolean ring (using AND as product and XOR as addition) and then converted easily to the diagrams shown. These `naive products' do define the product of any two vectors in $V$ but one should note that they do not define it on nondecomposable (`entangled') vectors such as $t\tens t+x\tens x$ since these are not in the image of the map $\tens:V\times V\to V\tens V$ (the image has 10 elements including 0). We would hardly worry about this in linear algebra since the product is linear but since we have not encoded such a property it is better to define the products more fully as maps $V\tens V\to V$ which we do in part (b). The products in part (a) factors through the maps in part (b) via the canonical map $\pi: V\times V\to V\tens V$ which as a diagram consists of 4 AND gates connecting up as shown in part (c). One can check with a little Boolean algebra that following this by the maps (b) gives the maps (a) so we can pull back to them, but the maps in (b) carry a little more information as explained. This language is obviously more tricky and for example associativity of the two ways to form the iterated product $V\tens V\tens V\to V$ ideally would be drawn in 4D with the input requiring a cube of wire-ends, one dimension for each tensor factor. In \cite{MaPac} we will describe differential structures and Riemannian geometry in this language and for small dimensions as outlined above.

\end{document}